\documentclass{amsart}
\usepackage{amsfonts}
\usepackage{upgreek}
\usepackage{eucal}
\usepackage{amssymb,epsfig,latexsym,amsthm}
\usepackage{graphicx,amsmath,amssymb,latexsym,psfrag,mathrsfs}
\usepackage{graphics,graphicx,color}
\usepackage{geometry}
\usepackage{diagbox}
\newtheorem{Theorem}{Theorem}[section]
\newtheorem{Proposition}[Theorem]{Proposition}

\newtheorem{Corollary}[Theorem]{Corollary}
\newtheorem{Lemma}[Theorem]{Lemma}
\newtheorem{Definition}[Theorem]{Definition}

\theoremstyle{definition}
\newtheorem{Remark}[Theorem]{Remark}

\newcommand{\dd}{\mathrm{d}}
\numberwithin{equation}{section}

\usepackage[
colorlinks=true,urlcolor=cyan,
citecolor=magenta,linkcolor=cyan,
linktocpage,pdfpagelabels,bookmarksnumbered,bookmarksopen,
pagebackref
]{hyperref}

\title{Normal approximation of random Gaussian neural networks}
\date{}
\makeatletter
\g@addto@macro{\endabstract}{\@setabstract}
\newcommand{\authorfootnotes}{\renewcommand\thefootnote{\@fnsymbol\c@footnote}}%
\newcommand{\authorsuperscript}{\@fnsymbol\c}%
\makeatother

\address[N.Apollonio, D. De Canditiis, G. Franzina, P. Stolfi, G. L. Torrisi]{\quad\newline
\indent Istituto per le Applicazioni del Calcolo ``M. Picone''
\newline\indent
Consiglio Nazionale delle Ricerche
\newline\indent 
Via dei Taurini, 19, 00185 Roma, Italy}


\begin{document}
\maketitle

\begin{center}
  \normalsize
  \authorfootnotes
  N. Apollonio\footnote{Corresponding author, \email{nicola.apollonio@cnr.it}}\textsuperscript{\ddag}, 
  D. De Canditiis\textsuperscript{\ddag},
  G. Franzina\textsuperscript{\ddag}, 
  P. Stolfi\textsuperscript{\ddag} and
  G.L. Torrisi\textsuperscript{\ddag} \par \bigskip

  \textsuperscript{\ddag}Istituto per le Applicazioni del Calcolo
  ``M. Picone'', CNR 

  \date{}
\end{center}
\begin{abstract}
In this paper we provide explicit upper bounds on 
some distances 
between the (law of the) output of
a random Gaussian neural network and (the law of) a random Gaussian vector. 
Our main results concern deep random Gaussian neural networks,  with a rather general activation function. The upper bounds show how the widths of the layers,  the activation function
and other architecture parameters affect the Gaussian approximation of the output.
Our techniques, relying on Stein's method and integration by parts formulas for the Gaussian law,  yield estimates on distances which are indeed integral probability metrics,  and include the convex distance.
This latter metric is defined by testing against indicator functions of measurable convex sets,
and so allow for accurate estimates of the probability that the output is localized
in some region of the space.  Such estimates have a significant interest both from a practitioner's and a theorist's perspective.
\end{abstract}

\tableofcontents

\section{Introduction}
This work is part of the literature studying random neural networks (NNs for short),  i.e., NNs whose biases and weights are random variables.  In the context
of modern deep learning, the interest in these types of networks is twofold: on the one hand, they  naturally constitute  a prior in a Bayesian approach, on the other hand they may
represent the initialization of gradient flows in empirical risk minimization,
see \cite{RYH} for a general reference on the subject.

Within the boundaries of this topic, many contributions in the literature have been handling the asymptotic Gaussianity of random NNs,  as the number of neurons in the hidden layers tends to infinity.  A seminal paper is \cite{Neal 1996},
where the output of a shallow (i.e.,  having a single hidden layer) random NN,  viewed as a stochastic process on the sphere,  is shown to converge to a Gaussian process, as the number of neurons in the hidden layer grows large.  From that
point onward,  many sophisticated results have been published for deep (i.e.,  having more than one hidden layer) random NNs. 
The starting point of our work is \cite{H},  where the output of a deep random NN, viewed as a random element on the space of continuous functions on a compact set,  is proved to converge to a Gaussian process, as the number of neurons in all the hidden layers tends to infinity.
Other related contributions are \cite{LBN18,Y,Y19}.


Recently,  the problem of the quantitative Gaussian approximation of the output of a random NN received a lot of attention. 
For instance,~\cite{EMS21}, exploiting Wasserstein distances,  provides quantitative versions of the results in
\cite{Neal 1996}
when the activation function is polynomial,  ReLu and hyperbolic tangent;
see \cite{BGRS, CMSV,  K22} for other important related contributions. In particular,  the reference \cite{BGRS} share with our work the idea to apply Stein's method for Gaussian approximation in the context of deep NN,  although in a different mathematical setting.

A significant achievement is provided in \cite{BT} where,  for the first time in the literature,  is given a quantitative proof of the Gaussian behavior of the output of a 
deep random Gaussian NN
(i.e.  a random NN whose biases and weights are Gaussian distributed)
with a Lipschitz continuous activation function; in \cite{BT} the distance from Gaussianity is measured
by means of 
the 2-Wasserstein metric,  which comes from the Monge-Kantorovich problem with quadratic cost. 
As far as shallow random Gaussian NNs with univariate output is concerned,  we mention the recent manuscript \cite{BFF23}
which provides quantitative bounds on the Kolmogorov,  the total variation and the 1-Wasserstein distances between the output and a Gaussian random variable,
when the activation function is sufficiently smooth and has a sub-polynomial growth.
A special mention deserves the independently written paper \cite{FHMNP} where Stein's method is used to obtain tight probabilistic bounds for various distances between the output (and its derivatives)
of a deep random Gaussian NN and a Gaussian random variable.  We refer the reader to Remarks \ref{re:18092023primo},  \ref{re:18092023secondo} and \ref{re:19092023} for comparisons between our results and the corresponding achiements in \cite{FHMNP}.

The main contribution of this paper concerns the Gaussian approximation of the output of  deep random Gaussian NNs in the convex and $1$-Wasserstein distances,  under mild assumptions on the activation function (which, differently from \cite{BT},  can be non-Lipschitz).  A specialization of these results clearly provides approximations for shallow random Gaussian NN with univariate output.  However, in this specific case we furnish direct proofs,  which (for various technical reasons) give the same rates under more general assumptions on the activation function. 

For shallow random Gaussian NNs with univariate output,
combining the Stein method for the Gaussian approximation
with the integration by parts formula for the Gaussian law,  we provide explicit bounds 
for the Kolmogorov,  the total variation and the $1$-Wasserstein
distances between the output and a Gaussian random variable,  under a minimal assumption on the activation function (see Theorem \ref{prop:shallow}).  Remarkably,
we obtain the same rate of convergence as in~\cite{BFF23},  as the number of neurons in the hidden layer grows large,  our constants being presumably better than the ones in \cite{BFF23} (see~Table~\ref{tab:shallow-fav-gian}). For
deep random Gaussian NNs,
the novelty of our results 
is that we measure the error in the Gaussian approximation of the output in terms of
the convex distance (see~Theorem~\ref{thm:convexdistance}) and of the
$1$-Wasserstein distance (see Theorem~\ref{thm:1Wasserstein}),  for a class of activation functions which includes the family of Lipschitz continuous functions (see Proposition \ref{prop:04052023primo}$(ii)$).  Remarkably,  for both the convex and the $1$-Wasserstein distances the rate of convergence that we obtain is of the same order as the one in~\cite{BT},  as the number of neurons in all the hidden layers tend to infinity.
The proofs of Theorems \ref{thm:convexdistance} and \ref{thm:1Wasserstein} are based on the Stein method for the multivariate Gaussian approximation and the integration by parts formula for the multivariate Gaussian law. The presence of more than one hidden layers complicates the derivations of the bounds, which rely on
a key estimate for the $L^2$-distance between the so-called collective observables and their limiting values (see Theorem \ref{le:fund_est}).  We emphasize that,
when considering the convex distance, 
an expedient tool is provided by a smoothing lemma that we borrow from~\cite{SY}.

It is well-known that localizing the output of a random NN, i.e.,  having a control over the probability that the output lies in a region of the space (belonging to a large class of measurable sets),
is of valuable interest for practitioners. 
From a theorist's point of view,  
the output distribution of a NN is often analytically untractable,  and computing the probability that the output belongs to some measurable set resorts in performing a \lq\lq heroic\rq\rq\,  mathematical integration (see e.g.\cite[p.~49]{RYH}).  Our 
Theorems \ref{prop:shallow}$(ii)$ and \ref{thm:convexdistance} offer some insights into the localization problem in a simple and efficient way for both the univariate output of a shallow random Gaussian NN and the output of a deep random Gaussian NN, respectively.  We refer the reader to Section \ref{sec:illustrations} for some numerical illustrations of this issue.

The paper is organized as follows.  In Section \ref{SEC: preliminary} we introduce our toolkit such as the integral probability metrics considered in the paper and some preliminaries on the Stein method.  In
Section~\ref{subsec:GNN},  first we introduce all the NNs considered in this work and then 
we give a brief overview
of the main results in \cite{BT} and \cite{BFF23}, comparing them with our achievements.
In Section~\ref{sec4} we give upper bounds for the Kolmogorov,  the total variation and the $1$-Wasserstein distances between the (univariate) output
of a shallow random Gaussian NN 
and a Gaussian random variable.
In
Section~\ref{sec:collective} we prove the aforementioned key estimate on the $L^2$-distance between the collective observables and their limiting values.  In Section \ref{sec:BOUNDS} we furnish explicit upper bounds on the convex and the $1$-Wasserstein distances between the output of a
deep random Gaussian NN and a Gaussian random vector.
Finally,  in Section \ref{sec:illustrations} we present some numerical illustrations concerning the above mentioned issue of the output localization.
\section{Preliminaries} \label{SEC: preliminary}

In the present section we introduce some notation and we recall
some results that will be of use throughout the paper.
\subsection{Distances between probability measures}
\label{sec:tools}
In this paper we consider various distances between probability measures on $\mathbb R^d$,  $d\in\mathbb N:=\{1,2,\ldots\}$: the total variation distance,  the convex distance,  the Komogorov distance, and the $p$-Wasserstein distances. Hereon, the symbol $\|\cdot\|_d$ denotes the Euclidean norm on $\mathbb R^d$. 

\begin{Definition}\label{df:TV}
The total variation distance between the laws of two $\mathbb R^d$-valued random vectors $\bold X$ and $\bold Y$,  written
$d_{TV}(\bold X,\bold Y)$, is given by
\[
d_{TV}(\bold X,\bold Y):=\sup_{B\in\mathcal{B}(\mathbb R^d)}|\mathbb{P}(\bold X\in B)-\mathbb{P}(\bold Y\in B)|,
\]
where $\mathcal{B}(\mathbb R^d)$ denotes the Borel $\sigma$-field on $\mathbb R^d$.
\end{Definition}

\begin{Definition}\label{df:C}
The convex distance between the laws of two $\mathbb R^d$-valued random vectors $\bold X$ and $\bold Y$, written
$d_c(\bold X,\bold Y)$, is given by
\[
d_c(\bold X,\bold Y):=\sup_{g\in\mathcal{I}_d}|\mathbb{E}[g(\bold X)]-\mathbb{E}[g(\bold Y)]|,
\]
where $\mathcal{I}_d$ denotes the collection of all indicator functions of measurable convex sets in $\mathbb{R}^d$.
\end{Definition}

\begin{Definition}\label{df:K}
The Kolmogorov distance between the laws of two $\mathbb R^d$-valued random vectors $\bold X=(X_1,\ldots,X_d)$ and $\bold Y=(Y_1,\ldots,Y_d)$,  written
$d_K(\bold X,\bold Y)$, is given by
\[
d_K(\bold X,\bold Y):=\sup_{\bold y=(y_1,\ldots,y_d)\in\mathbb R^d}|\mathbb{P}(X_1\leq y_1,\ldots,X_d\leq y_d)-\mathbb{P}(Y_1\leq y_1,\ldots,Y_d\leq y_d)|.
\]
\end{Definition}

\begin{Definition}\label{def:2Wasserstein}
For $p\in[1,+\infty)$, the $p$-Wasserstein distance between the laws of two $\mathbb R^d$-valued random vectors $\bold X$ and $\bold Y$, 
written
$d_{W_p}(\bold X,\bold Y)$, is given by
\[
d_{W_p}(\bold X,\bold Y):=\inf_{(\bold U,\bold V)\in C_{(\bold X,\bold Y)}}\mathbb{E}[\|\bold U-\bold V\|_d^p]^{1/p},
\]
where $C_{(\bold X,\bold Y)}$ is the family of all the couplings of $\bold X$ and $\bold Y$, i.e.,  the family of all random vectors $(\bold U,\bold V)$ such that $\bold U$ is distributed as $\bold X$ and $\bold V$ is distributed as $\bold Y$.
\end{Definition}

Clearly,  by Jensen's inequality $d_{W_1}\leq d_{W_2}$, and it follows directly by the definitions that $d_K\leq d_c\leq d_{TV}$.  Furthermore, 
for all $s=TV,c,K, W_p$,
if $d_s(\bold{Y}_n,\bold Y)\to 0$,  as $n\to+\infty$, where $\bold{Y}_n$,  $n\in\mathbb N$,  and $\bold Y$ are random vectors with values in $\mathbb R^d$,  then $\bold Y_n$ converges in law to $\bold Y$,  as $n\to+\infty$
(see e.g.  \cite{NPbook,V}).

In view of the Kantorovich-Rubinstein duality
(see Theorem~5.10 and Eq.\ (5.11) in \cite{V}),  the $1$-Wasserstein distance between the laws of two $\mathbb R^d$-valued random vectors $\bold X$ and $\bold Y$ such that $\max\{\mathbb E\|\bold X\|_d,\mathbb E\|\bold Y\|_d\}<\infty$,  satisfies the relation
\begin{equation}
\label{earthmov}
d_{W_1}(\bold X,\bold Y)=\sup_{g\in\mathcal{L}_d(1)}|\mathbb{E}[g(\bold X)]-\mathbb{E}[g(\bold Y)]|,
\end{equation}
where $\mathcal{L}_d(1)$ indicates the collection of all functions $g:\mathbb R^d\to\mathbb R$ which are Lipschitz continuous with Lipschitz constant less than or equal to $1$.
Since $d_c$ is defined by testing against 
indicator functions of Borel convex sets rather
than arbitrary Borel sets, the convex distance can be expected
to be estimated easier than the total variation distance; moreover,  the convex distance looks more flexible than the Kolmogorov distance,  for example it enjoys a number of invariance properties not satisfied by $d_K$,  see~\cite{Benktus}. 

As for the relation between the convex distance
and the optimal transport metric $d_{W_1}$,  it turns out that
the convex distance to a fixed centered Gaussian law is bounded
from above by a multiple of the square root of the $1$-Wasserstein distance.  More precisely, one has the following Proposition~\ref{prop:dclessdw}, which is proved in \cite{NPX}.

Here and henceforth,
we denote by $\bold{N}_\Sigma=(N_1,\ldots,N_d)$,  $d\in\mathbb N$,  a centered Gaussian vector with invertible covariance matrix $\Sigma=(\Sigma_{ij})_{1\leq i,j\leq d}$.  

\begin{Proposition}\label{prop:dclessdw}
For any $d$-dimensional random vector $\bold Y$,  we have
\[
d_c(\bold Y,\bold{N}_\Sigma)\leq 2\sqrt 2\Gamma(\Sigma)^{1/2}d_{W_1}(\bold Y,\bold{N}_\Sigma)^{1/2},
\]
where $\Gamma(\Sigma)$ is the constant defined by
\begin{equation}\label{eq:isoperimetric}
\Gamma(\Sigma):=\sup_{Q,\,\epsilon>0}\frac{\mathbb{P}(\bold{N}_\Sigma\in Q^\epsilon)-\mathbb{P}(\bold N_\Sigma\in Q)}{\epsilon},
\end{equation}
where $Q$ ranges over all the Borel measurable convex subsets of $\mathbb R^d$,  and $Q^\epsilon$ denotes the set of all elements of $\mathbb R^d$ whose Euclidean distance from $Q$ does not exceed $\epsilon$.
\end{Proposition}

$\Gamma(\Sigma)$ is in fact an isoperimetric constant.  Indeed,  if $Q$ is a bounded convex Borel set, then it is easily seen that
\[
\sup_{\epsilon>0}\frac{\mathbb{P}(\bold{N}_\Sigma\in Q^\epsilon)-\mathbb{P}(\bold N_\Sigma\in Q)}{\epsilon}
= \int_{\partial Q} h_\Sigma\left(\mathbf x,\boldsymbol{\nu}^Q_{\mathbf x}\right)d\,\mathscr H^{d-1}\,,
\]
where $\mathscr H^{d-1}$ is the $(d-1)$-dimensional Hausdorff
measure, for $\mathscr H^{d-1}$-a.e.\ boundary point $\mathbf x$
we are denoting by $\boldsymbol{\nu}_{\mathbf x}^Q$ the
outward unit normal to $\partial Q$ at $\mathbf x$, and we
have set
\[
	h_\Sigma(\mathbf x,\boldsymbol{\nu}) := (2\pi)^{-d/2} 
	e^{-\frac{|\mathbf x|^2}{2}}\|\Sigma^{-\frac12} \boldsymbol{\nu}_{\mathbf x}^Q\|_d\,,\qquad \text{for all $(\mathbf x\mathbin,\boldsymbol{\nu})
	\in \mathbb R^d\times\mathbb S^{d-1}$,}
\]
where $\mathbb{S}^{d-1}$ is the unit sphere of $\mathbb R^d$.
Thus,  $\Gamma(\Sigma)$ is
the worst (largest) possible anisotropic Gaussian perimeter for a 
convex body in $\mathbb R^d$.

We refer the reader to~\cite{naz04} for the following bound
\begin{equation}
\label{worstgauss}
e^{-\frac54} d^{1/4} \le \Gamma(\Sigma) \le (2\pi)^{-\frac14}d^{1/4}.
\end{equation}

\subsection{The one-dimensional Stein equation}

Throughout this paper,  we denote by $\mathcal{N}(\mu,\eta)$ the one-dimensional Gaussian law with mean $\mu\in\mathbb R$ and variance $\eta>0$,  and let
$Z\sim\mathcal{N}(0,1)$.  

The celebrated Stein equation for the one-dimensional Normal approximation \cite{CS} is given by
\begin{equation}\label{eq:stein1d}
g(w)-\mathbb Eg(Z)=f_g'(w)-wf_g(w),
\end{equation}
where $g:\mathbb R\to\mathbb R$ is a measurable function such that $\mathbb E |g(Z)|<\infty$ and $f_g:\mathbb R\to\mathbb R$ is unknown.  The following lemma holds, see e.g.  Proposition 3.2.2,  Theorem 3.3.1,  Theorem 3.4.2 and Proposition 3.5.1 in \cite{NPbook}. See \cite{CGQ} for an introduction on the Stein method.

Hereon, for a Lipschitz continuous function $g:\mathbb R^d\to\mathbb R$ we denote by $\mathrm{Lip}(g)$ the Lipschitz constant of $g$,  and for a function $f:\mathbb R^d\to\mathbb R$ we denote by $\|f\|_\infty$ the supremum norm of $f$.

\begin{Lemma}\label{le:CGQ}
The following claims hold:\\
\noindent$(i)$ For any $y\in\mathbb R$,  the Stein equation \eqref{eq:stein1d} with $g(w):=\bold{1}_{(-\infty,y]}(w)$ has a unique solution $f_g$ 
and 
$\|f_g'\|_\infty\leq 1$. \\
\noindent$(ii)$ Let $g:\mathbb R\to [0,1]$ be a measurable function. Then there exists a unique solution $f_g$ of the Stein equation \eqref{eq:stein1d} and
$\|f'_g\|_\infty\leq 2.$\\
\noindent$(iii)$ Let $g:\mathbb R\to\mathbb R$ be a Lipschitz continuous function.  Then there exists a unique solution $f_g$ of the Stein equation \eqref{eq:stein1d} and
$\|f'_g\|_\infty\leq\mathrm{Lip}(g)\sqrt{2/\pi}.$
\end{Lemma}

\subsection{The multidimensional Stein equation}

Throughout this paper,  given a sufficiently smooth function $f:\mathbb R^d\to\mathbb R$,  we define 
\[
\partial_{i_1 i_1\ldots i_n}^{n} f(x_1,\ldots,x_d):=\frac{\partial^n f}{\partial x_{i_1}\ldots\partial x_{i_n}}(x_1,\ldots,x_d).
\]

Let $\mathcal{M}_{d\times d}(\mathbb R)$,  $d\in\mathbb N$,  be the set of
$d\times d$ real matrices.  For a function $f\in C^2(\mathbb R^d)$,  we denote by  $\mathrm{Hess}\,f(\bold y)\in\mathcal{M}_{d\times d}(\mathbb R)$ the Hessian matrix of $f$ at $\bold y\in\mathbb R^d$ and by $\|\cdot\|_{op}$ the operator norm on $\mathcal{M}_{d\times d}(\mathbb R)$,  i.e.,  for any $\Gamma\in\mathcal{M}_{d\times d}(\mathbb R)$,  $\|\Gamma\|_{\mathrm{op}}:=\sup_{\bold y:\,\,\|\bold y\|_d=1}\|\Gamma\bold y\|_d$.  We consider the Hilbert-Schmidt inner product and the Hilbert-Schmidt norm on $\mathcal{M}_{d\times d}(\mathbb R)$,  which are defined,  respectively,  by
\[
\langle\Gamma,\Psi\rangle_{H.S.}:=\mathrm{Tr}(\Gamma \Psi^\top)=\sum_{i,j=1}^{d}\Gamma_{ij}\Psi_{ij}\quad\text{and}\quad\|\Gamma\|_{H.S.}=\sqrt{\langle\Gamma,\Gamma\rangle_{H.S.}}
\]
for every pair of matrices $\Gamma=(\Gamma_{ij})_{1\leq i,j\leq j}$ and $\Psi=(\Psi_{ij})_{1\leq i,j\leq d}$, where the symbols $\mathrm{Tr}(\Gamma)$ and $\Gamma^\top$ denote, respectively,  the trace and the transpose of the matrix $\Gamma$. 

The Stein equation for multivariate Normal approximation is defined as

\begin{equation}\label{eq:SteinNew}
g(\bold y)-\mathbb{E}[g(\bold N_\Sigma)]=\langle\bold y,\nabla f_g(\bold y)\rangle_d-\langle\Sigma,\mathrm{Hess}\,f_g(\bold y)\rangle_{H.S.},\quad\bold y\in\mathbb R^d
\end{equation}
where $g:\mathbb R^d\to\mathbb R$ is given and $f_g$ is unknown.  Throughout this paper the symbol $\langle\cdot,\cdot\rangle_d$ denotes the inner product in $\mathbb R^d$.

The following lemmas provide solutions to Stein's equation \eqref{eq:SteinNew},  under different assumptions on $g$.  

\begin{Lemma}\label{le:Stein}
Let $g\in\mathcal{L}_d(1)$. Then the function
\[
f_g(\bold y):=\int_0^\infty\mathbb{E}[g(\bold{N}_\Sigma)-g(\mathrm{e}^{-t}\bold y+\sqrt{1-\mathrm{e}^{-2t}}\bold N_\Sigma)]\,\dd t,\quad \bold y\in\mathbb R^d
\]
is such that: $f_g\in C^2(\mathbb R^d)$,  $f_g$ solves \eqref{eq:SteinNew},  $f_g$ satisfies
\begin{equation}\label{eq:18052023secondo}
\|\partial_i f_g\|_\infty\leq 1,\quad\text{for any $i=1,\ldots,d$}
\end{equation}
and
\[
\sup_{\bold{y}\in\mathbb R^d}\|\mathrm{Hess}f_g(\bold y)\|_{H.S.}\leq\sqrt{d}\|\Sigma^{-1}\|_{op}\|\Sigma\|_{op}^{1/2}.
\]
\end{Lemma}

\begin{Lemma}\label{le:Stein2}
For $g:\mathbb{R}^d\to\mathbb{R}$ measurable and bounded,  define the smoothed function
\begin{equation}
\label{eq:10_01_2022}
g_t(\bold y):=\mathbb{E}[g(\sqrt{t}\bold N_{\Sigma}+\sqrt{1-t}\bold y)],\quad\text{$\bold y\in\mathbb{R}^d$, }
\end{equation}
where $t\in (0,1)$ is a smoothing parameter.
Then:\\ 
\noindent$(i)$ For any $t\in (0,1)$,  the function
\[
f_{t,g}(\bold y):=\frac{1}{2}\int_t^1\frac{1}{1-s}\mathbb{E}[g(\sqrt s \bold N_\Sigma+\sqrt{1-s}\bold y)-g(\bold N_\Sigma)]\,\dd s,\quad\bold y\in\mathbb{R}^d
\]
is such that: $f_{t,g}\in C^2(\mathbb R^d)$,  $f_{t,g}$ solves \eqref{eq:SteinNew} with $g_t$ in place of $g$, and
\begin{equation}\label{eq:17052023sesto}
\|\partial_i f_{t,g}\|_\infty\leq\|g\|_\infty\sqrt{\frac{1-t}{t}}\sum_{\ell,j=1}^{d}(\Sigma^{-1/2})_{\ell j}(\Sigma^{-1/2})_{\ell i}\sqrt{\Sigma_{jj}},
\quad\text{for any $i=1,\ldots,d$.}
\end{equation}
\noindent$(ii)$ For any $d$-dimensional random vector $\bold Y$ it holds
\[
\sup_{g\in\mathcal{I}_d}\mathbb E\|\mathrm{Hess}(f_{t,g}(\bold Y))\|_{H.S.}^2
\leq\|\Sigma^{-1}\|_{op}^{2}(d^2(\log t)^2 d_c(\bold Y,\bold N_{\Sigma})+530 d^{17/6}),\quad\text{for any $t\in (0,1)$.}
\]
\end{Lemma}

See Proposition 4.3.2 in \cite{NPbook} for Lemma \ref{le:Stein};
in particular in 
 \cite{NPbook} it is noticed that
\[
\partial_i f_g(\bold y)=-\int_0^{\infty}\mathrm{e}^{-t}\mathbb{E}\left[\partial_i g(\mathrm{e}^{-t}\bold y+\sqrt{1-\mathrm{e}^{-2t}}\bold{N}_{\Sigma})\right]\mathrm{d}t,
\quad\text{$i=1,\ldots,d$}
\]
which, combined with the fact that $g\in\mathcal{L}_d(1)$,  gives  
the bound \eqref{eq:18052023secondo}; see
\cite{SY} p. 12 and Proposition 2.3 for Lemma \ref{le:Stein2},  and Lemma 3.6$(ii)$ in \cite{To} for the bound \eqref{eq:17052023sesto}.

\subsection{The smoothing lemma and the integration by parts formula for Gaussian random vectors}

We state a remarkable smoothing lemma for the convex distance proved in \cite{SY}, see Lemma 2.2 therein.  It plays a crucial role in the proof of the Normal approximation of the output of a deep random Gaussian NN in the metric $d_c$, see Theorem \ref{thm:convexdistance}.

\begin{Lemma}\label{le:smoothingdconvex}
Let $\bold Y$ be a $d$-dimensional random vector.  Then,  for any $t\in (0,1)$,  
\[
d_c(\bold Y,\bold N_{\Sigma})\leq\frac{4}{3}\sup_{g\in\mathcal{I}_d}|\mathbb{E}[g_t(\bold Y)-g_t(\bold N_\Sigma)]|+\frac{20 d}{\sqrt 2}\frac{\sqrt t}{1-t},
\]
where $g_t$ is defined by \eqref{eq:10_01_2022}.
\end{Lemma} 

We recall the Gaussian integration by parts formula (we refer the reader 
to Exercise 3.1.4 in \cite{NPbook} for the Part $(i)$ of Lemma \ref{le:Gaussibp}
and to  Exercise 3.1.5 in \cite{NPbook} for the Part $(ii)$ of Lemma \ref{le:Gaussibp}).
\begin{Lemma}\label{le:Gaussibp} 
The following claims hold:\\
\noindent$(i)$ $N\sim\mathcal{N}(\mu,\eta)$ if and only if,  for any differentiable function $g:\mathbb R\to\mathbb R$ such that $\mathbb E |g'(N)|<\infty$,  we have $\mathbb E(N-\mu)g(N)=\eta\mathbb Eg'(N)$.\\
\noindent$(ii)$ Let $g\in C^1(\mathbb R^d)$ with bounded first partial derivatives. 
Then
\[
\mathbb{E}[N_i g(\bold N_\Sigma)]=\sum_{j=1}^{d}\Sigma_{ij}\mathbb{E}[\partial_j g(\bold N_\Sigma)],
\quad\text{for any $i=1,\ldots,d$.}
\]

\end{Lemma}
%

\section{Random neural networks}\label{subsec:GNN}

We let $L\in\mathbb N$,  we take
 $L+2$ positive integers $n_0,\ldots,n_{L+1}\in\mathbb N$, and we fix a function 
$\sigma:\mathbb R\to\mathbb R$.  
A fully connected NN of depth $L$ with input dimension $n_0$,  output dimension $n_{L+1}$,
hidden layer widths $n_1,\ldots,n_L$ and non-linearity $\sigma$ is a mapping
\[
\bold x:=(x_1,\ldots,x_{n_0})\in\mathbb R^{n_0}\mapsto\bold{z}^{(L+1)}(\bold x)=(z_1^{(L+1)}(\bold x),\ldots,z_{n_{L+1}}^{(L+1)}(\bold x)))\in\mathbb{R}^{n_{L+1}}
\] 
that is defined by a recursive relation of the form
\begin{align}
z_i^{(1)}(\bold x)&=b_i^{(1)}+\sum_{j=1}^{n_0}W_{ij}^{(1)}x_j,\quad\text{$i=1,\ldots,n_{1}$}\nonumber\\
z_i^{(\ell)}(\bold x)&=b_i^{(\ell)}+\sum_{j=1}^{n_{\ell-1}}W_{ij}^{(\ell)}\sigma(z_j^{(\ell-1)}(\bold x)),\quad\text{$i=1,\ldots,n_{\ell}$,  $\ell=2,\ldots,L+1$}\nonumber
\end{align}
where the parameters $b_i^{(\ell)}\in\mathbb R$ and $W_{ij}^{(\ell)}\in\mathbb R$ are called network biases and weights,  respectively. The quantities $L$ and $n_0,\ldots,n_{L+1}$ constitute the so-called network architecture.  The function $\sigma$ is usually called activation function. 
NNs of this kind will be denoted by
\[
\mathrm{NN}(L,n_0,n_{L+1},\bold{n}_L,\sigma,\bold x,\bold b,\bold W)\,,
\]
 where $\bold{n}_L:=(n_1,\ldots,n_{L})$,
$\bold{b}:=(b_{i}^{(\ell)})$ and $\bold{W}:=(W_{ij}^{(\ell)})$.

We say that the neural network $\mathrm{NN}(L,n_0,n_{L+1},\bold n_L,\sigma,\bold x,\bold b,\bold W)$ is a (fully connected and) deep random Gaussian neural network, denoted by 
\[
\mathrm{GNN}(L,n_0,n_{L+1},\bold n_L,\sigma,\bold x,\bold b,\bold W)\,,
\]
if 
$\sigma:\mathbb R\to\mathbb R$ is measurable and
$b_i^{(\ell)},W_{ij}^{(\ell)}$, $i=1,\ldots,n_\ell$, $j=1,\ldots,n_{\ell-1}$, $\ell=1,\ldots,L+1$,  are independent random variables with 
\[
b_i^{(\ell)}\sim\mathcal{N}(0,C_b)\quad\text{and}\quad W_{ij}^{(\ell)}\sim\mathcal{N}(0,C_W/n_{\ell-1}),\quad\text{$\ell=1,\ldots,L+1$}
\]
for positive constants $C_b,C_W>0$. 

NNs of depth $L=1$
are called shallow NNs.
We shall denote shallow NNs
(respectively,  shallow random Gaussian NNs) by 
$
\mathrm{NN}(1,n_0,n_{2},n_1,\sigma,\bold x,\bold b,\bold W)$ 
(respectively,  by $\mathrm{GNN}(1,n_0,n_2,n_1,\sigma,\bold x,\bold b,\bold W)$). 

Throughout this paper we will consider also NNs with univariate output, i.e., NNs with
$n_{L+1}=1$. 

Consider a deep random Gaussian neural network $\mathrm{GNN}(L,n_0,n_{L+1},\bold n_L,\sigma,\bold x,\bold b,\bold W)$.
It turns out that the random variables $z_i^{(1)}=z_i^{(1)}(\bold x)$, $i=1,\ldots,n_1$,  are independent and identically distributed with
\[
z_i^{(1)}\sim\mathcal{N}\left(0,C_b+\frac{C_W}{n_0}\sum_{j=1}^{n_0}x_j^2\right).
\]
For $\ell=1,\ldots,L$,  let $\mathcal{F}_\ell$ be the $\sigma$-field generated by the random variables
\[
\{b_i^{(h)},W_{ij}^{(h)},\quad i=1,\ldots,n_h,\, j=1,\ldots,n_{h-1},\,h=1,\ldots,\ell\}.
\]
By construction,  for any fixed $\ell\in\{2,\ldots,L+1\}$,
given $\mathcal{F}_{\ell-1}$,  the random variables $z_i^{(\ell)}=z_i^{(\ell)}(\bold x)$, $i=1,\ldots,n_{\ell}$,  are independent and Gaussian (as linear combination of independent Gaussian random variables).  
A straightforward computation yields
\[
\mathbb{E}[z_i^{(\ell)}\,|\,\mathcal{F}_{\ell-1}]=0,\quad\text{$i=1,\ldots,n_{\ell}$}
\]
and
\begin{equation}\label{eq:17052023quarto}
\mathbb{E}[|z_i^{(\ell)}|^2\,|\,\mathcal{F}_{\ell-1}]=C_b+\frac{C_W}{n_{\ell-1}}\sum_{j=1}^{n_{\ell-1}}|\sigma(z_j^{(\ell-1)})|^2,\quad\text{$i=1,\ldots,n_{\ell}$.}
\end{equation}

Setting $\bold{n}_{\ell}=(n_1,\ldots,n_{\ell})$,  $\ell=1,\ldots,L$,  we define the quantities
\[
\mathcal{O}_{\bold{n}_\ell}^{(\ell)}:=\frac{1}{n_{\ell}}\sum_{j=1}^{n_\ell}\sigma(z_j^{(\ell)})^2,\quad\text{$\ell=1,\ldots,L$}
\]
and 
\begin{equation}\label{eq:20062023sesto}
\mathcal{O}^{(\ell)}:=\mathbb{E}\left[\sigma\left(Z\sqrt{C_b+C_W\mathcal{O}^{(\ell-1)}}\right)^2\right],\quad\text{$\ell=1,\ldots,L$}
\end{equation}
where
\begin{equation}\label{eq:20062023quinto}
\mathcal{O}^{(0)}:=\frac{1}{n_0}\sum_{j=1}^{n_0}x_j^2\quad\text{and}\quad Z\sim\mathcal{N}(0,1).
\end{equation}
In the literature, the random variable $\mathcal{O}_{\bold{n}_\ell}^{(\ell)}$ is often referred to  as collective  observable at layer $\ell$.

\subsection{Some related literature}
Consider the output
$
\bold{z}^{(L+1)}:=(z_1^{(L+1)},\ldots,z_{n_{L+1}}^{(L+1)})
$
of a deep random Gaussian neural network $\mathrm{GNN}(L,n_0,n_{L+1},\bold{n}_{L},\sigma;\bold x;\bold b,\bold W)$ and let
\[
\bold{z}:=(z_1,\ldots,z_{n_{L+1}})
\]
be a centered $n_{L+1}$-dimensional Gaussian random vector with (invertible) covariance matrix
\begin{equation}\label{eq:09052023primo}
\Sigma_{n_{L+1}}:=(C_b+C_W\mathcal{O}^{(L)})\mathrm{Id}_{n_{L+1}},
\end{equation}
where $\mathrm{Id}_{n_{L+1}}$ is the identity matrix of $\mathcal{M}_{n_{L+1}\times n_{L+1}}(\mathbb R)$.  It follows from  Theorem 1.2 in \cite{H} (which indeed, more generally,  establishes a functional weak convergence) that,  if $\sigma$ is continuous and polynomially bounded,  then
\begin{equation}\label{eq:hanin}
\bold{z}^{(L+1)}\to \bold z\quad\text{in law,  as $\min\{n_1,\ldots,n_L\}\to+\infty$.}
\end{equation}

The following result for shallow random Gaussian NNs is proved in \cite{BFF23}, see Theorem 3.2 therein.

\begin{Theorem}\label{thm:BFF23}
Let $\mathrm{GNN}(1,n_0,1,n_1,\sigma,\bold x,\bold b,\bold W)$ be a shallow random Gaussian NN with univariate output.  If
\begin{equation}\label{eq:favaro}
\text{$\sigma\in C^2(\mathbb R)$ and $\max\{|\sigma(x)|,|\sigma'(x)|,|\sigma''(x)|\}\leq r_1+r_2|x|^\gamma$, $x\in\mathbb R$}
\end{equation}
for some $r_1,r_2,\gamma\geq 0$,  then
\begin{align}
d_s(z^{(2)},z)&\leq c_s\sqrt{C_b+C_W\mathcal{O}^{(0)}+(C_b+C_W\mathcal{O}^{(0)})^2(2+\sqrt{3(1+2(C_b+C_W\mathcal{O}^{(0)})+3(C_b+C_W\mathcal{O}^{(0)})^2)})}\nonumber\\
&\qquad\qquad\qquad\qquad\qquad\qquad
\times\|r_1+r_2|Z\sqrt{C_b+C_W\mathcal{O}^{(0)}}|^\gamma\|_{L^4}^2\frac{1}{\sqrt{n_1}},
\end{align}
where 
$s=TV,K,W_1$ and
\[
c_{TV}:=\frac{4}{C_b+C_W\mathcal{O}^{(1)}},\quad c_K;=\frac{2}{C_b+C_W\mathcal{O}^{(1)}},\quad c_{W_1}:=\frac{1}{\sqrt{C_b+C_W\mathcal{O}^{(1)}}}\sqrt{8/\pi}.
\]
\end{Theorem}

In Section \ref{sec4},  we will give bounds on the quantities $d_s(z^{(2)},z)$,  $s=TV,K,W_1$,  of order $1/\sqrt{n_1}$, as $n_1\to\infty$,  under a minimal assumption on the activation function (note that Condition \eqref{eq:favaro} excludes the important case of the ReLu function, i.e.,  $\sigma(x):=x\bold{1}\{x\geq 0\}$),  see Theorem \ref{prop:shallow}.  
In Section \ref{sec:BOUNDS} we will give two general bounds on $d_{c}(\bold{z}^{(L+1)},\bold{z})$ and $d_{W_1}(\bold{z}^{(L+1)},\bold{z})$ for deep random Gaussian NNs, see Theorems \ref{thm:convexdistance} and \ref{thm:1Wasserstein}.
When specialized to shallow random Gaussian neural networks $\mathrm{GNN}(1,n_0,n_2,n_1,\sigma,\bold x,\bold b,\bold W)$,  they provide computable bounds respectively on $d_{c}(\bold{z}^{(2)},\bold{z})$ and $d_{W_1}(\bold{z}^{(2)},\bold{z})$ of order $1/\sqrt n_1$,  as $n_1\to\infty$,  
see Theorem 3.3 in \cite{BFF23} for a related result.

The first result in literature which quantifies the convergence in distribution \eqref{eq:hanin} with $L\geq 2$ is given in \cite{BT},  where the following theorem has been proved.  

\begin{Theorem}\label{thm:trevisan}
Let $\mathrm{GNN}(L,n_0,n_{L+1},\bold{n}_L,\sigma,\bold x,\bold b,\bold W)$ be a deep random Gaussian NN,  and suppose that the activation function $\sigma$ is Lipschitz continuous. Then
\[
d_{W_2}(\bold{z}^{(L+1)},\bold{z})\leq\sqrt{n_{L+1}}\sum_{i=1}^{L}\frac{C^{(i+1)}[\mathrm{Lip}(\sigma)\sqrt{C_W}]^{L-i}}{\sqrt{n_i}},
\]
where, for any $i=1,\ldots,L$,  $C^{(i+1)}$ are explicitly known positive constants, depending upon $\sigma$, $\bold x$, $C_b$ and $C_W$.
\end{Theorem}

The next corollary is an immediate consequence of Theorem \ref{thm:trevisan},  the 
fact that $d_{W_1}\leq d_{W_2}$,  Proposition \ref{prop:dclessdw} and \eqref{worstgauss}.

\begin{Corollary}\label{cor:19062023}
Let the asssumptions and notation of Theorem \ref{thm:trevisan} prevail. Then
\begin{equation}\label{eq:dWalternative}
d_{W_1}(\bold{z}^{(L+1)},\bold{z})\leq\sqrt{n_{L+1}}\sum_{i=1}^{L}\frac{C^{(i+1)}[\mathrm{Lip}(\sigma)\sqrt{C_W}]^{L-i}}{\sqrt{n_i}}
\end{equation}
and
\begin{equation}\label{eq:dcsub}
d_{c}(\bold{z}^{(L+1)},\bold{z})\leq 2^{\frac{11}{8}}
	\pi^{-\frac18} n_{L+1}^{\frac38}
\left(\sum_{i=1}^{L}\frac{C^{(i+1)}[\mathrm{Lip}(\sigma)\sqrt{C_W}]^{L-i}}{\sqrt{n_i}}\right)^{1/2}.
\end{equation}
\end{Corollary}

Theorems~\ref{thm:convexdistance} and~\ref{thm:1Wasserstein} will provide, under more general assumptions on $\sigma$ (see Proposition \ref{prop:04052023primo}$(ii)$) 
explicit bounds respectively on $d_{c}(\bold{z}^{(L+1)},\bold{z})$ and $d_{W_1}(\bold{z}^{(L+1)},\bold{z})$ which are of the same order of the bound in \eqref{eq:dWalternative},  as $n_1,\ldots,n_L\to\infty$.  In particular,  the bound on the convex distance of Theorem~\ref{thm:convexdistance}
considerably improves the one in \eqref{eq:dcsub}.  

\begin{Remark}
Although not strictly related to our results,  the pioneering papers \cite{EMS21} and \cite{Neal 1996} deserve a special mention.  In \cite{Neal 1996} the author considers a random shallow $NN(1,n_0,1,n_1,\sigma,\bold x,\bold 0,\bold W)$ with univariate output,  $b_i^{(1)}:=0$,  for all $i=1,\ldots,n_1$,  $b_1^{(2)}:=0$,
$W_{ij}^{(1)}$, $i=1,\ldots,n_1$, $j=1,\ldots,n_0$,  independent with law $\mathcal{N}(0,1)$ and $W_{1j}^{(2)}$,  $j=1,\ldots,n_1$,  independent, identically distributed with law $\mathbb{P}(W_{1j}^{(2)}=\pm 1/\sqrt{n_1})=1/2$ and independent of the random variables $W_{ij}^{(1)}$.  
It is proved in \cite{Neal 1996} that there exists a Gaussian process $G$ on $\mathbb{S}^{n_0-1}$ (the unit sphere in $\mathbb{R}^{n_0}$) such that the process $\{z^{(2)}(\bold x)\}_{\bold x\in\mathbb{S}^{n_0-1}}$ converges in distribution to $G$, as $n_1\to\infty$.  Quantitative versions of this result (for various Wasserstein metrics and some specific choices of $\sigma$) are provided in \cite{EMS21}.
\end{Remark}

\section{Normal approximation of shallow random Gaussian NNs with univariate output}\label{sec4}

The following theorem holds.

\begin{Theorem}\label{prop:shallow}
Let $\mathrm{GNN}(1,n_0,1,n_1,\sigma,\bold x,\bold b,\bold W)$ be a shallow random Gaussian NN with univariate output,  and assume that the activation function $\sigma$ is such that
\begin{equation}\label{eq:03072023}
0<\mathbb{V}\mathrm{ar}(\sigma(Z\sqrt{C_b+C_W\mathcal{O}^{(0)}})^2)<\infty.
\end{equation}
Then:\\
\noindent$(i)$ 
\[
d_K(z^{(2)},z)\leq\frac{C_W\sqrt{\mathbb{V}\mathrm{ar}(\sigma(Z\sqrt{C_b+C_W\mathcal{O}^{(0)}})^2)}}{C_b+C_W\mathbb{E}\sigma(Z\sqrt{C_b+C_W\mathcal{O}^{(0)}})^2}\frac{1}{\sqrt n_1}.
\]
\noindent$(ii)$ 
\[
d_{TV}(z^{(2)},z)\leq2\frac{C_W\sqrt{\mathbb{V}\mathrm{ar}(\sigma(Z\sqrt{C_b+C_W\mathcal{O}^{(0)}})^2)}}{C_b+C_W\mathbb{E}\sigma(Z\sqrt{C_b+C_W\mathcal{O}^{(0)}})^2}\frac{1}{\sqrt n_1}.
\]
\noindent$(iii)$
\[
d_{W_1}(z^{(2)},z)\leq\sqrt{2/\pi}\frac{C_W\sqrt{\mathbb{V}\mathrm{ar}(\sigma(Z\sqrt{C_b+C_W\mathcal{O}^{(0)}})^2)}}{\sqrt{C_b+C_W\mathbb{E}\sigma(Z\sqrt{C_b+C_W\mathcal{O}^{(0)}})^2}}\frac{1}{\sqrt n_1}.
\]

Note that the bound on the Kolmogorov distance is better than the one which can be obtained using the relation $d_K\leq d_{TV}$. 
\end{Theorem}

\begin{Remark}\label{re:18092023primo} 
Remarkably, Theorem 3.3 in \cite{FHMNP} shows that if $\mathrm{GNN}(1,n_0,1,n_1,\sigma,\bold x,\bold b,\bold W)$ is a shallow random Gaussian NN with univariate output,  and the activation function $\sigma$
is polynomially bounded to order $r\geq 1$ (see Definition 2.1 in \cite{FHMNP}),  then there exist two constants $C,C_0>0$ such that
\[
\frac{C_0}{n_1}\leq\max\{d_{W_1}(z^{(2)},z),d_{TV}(z^{(2)},z)\}\leq\frac{C}{n_1}.
\]
Although this inequality shows the optimality of the rate $1/n_1$, since the constants are not provided in closed form,  it can not be directly used for the purpose of output localization (see Section \ref{sec:illustrations}). Moreover,  the assumption \eqref{eq:03072023} on $\sigma$ does not require any regularity of the activation function.
\end{Remark}

\noindent{\it Proof.} We prove the three bounds $(i)$, $(ii)$, and $(iii)$ separately,
by the Stein method.\\
\noindent{\it Proof\,\,of\,\,Part\,\,$(i)$.} Set 
\begin{equation}\label{eq:030720231}
\nu:=\sqrt{C_b+C_W\mathcal{O}^{(1)}}, 
\end{equation}
and consider the Stein equation \eqref{eq:stein1d}
with 
\[
g(w):=\bold{1}_{(-\infty,y]}(\nu w)\,.
\]  Let $f_g$ be the unique solution of the Stein equation (see Lemma \ref{le:CGQ}$(i)$). Then,  for any $y\in\mathbb R$, 
\[
\bold{1}\{z^{(2)}\leq y\}-\mathbb{P}(z\leq y)=f'_{g}(z^{(2)}/\nu)-(z^{(2)}/\nu)f_{g}(z^{(2)}/\nu).
\]
Taking the expectation, we have
\[
\mathbb{P}(z^{(2)}\leq y)-\mathbb{P}(z\leq y)=\mathbb E[f'_{g}(z^{(2)}/\nu)-(z^{(2)}/\nu)f_{y/\nu}(z^{(2)}/\nu)].
\]
By Lemma \ref{le:Gaussibp}$(i)$ we have
\begin{align}
\mathbb{E}[(z^{(2)}/\nu)f_{g}(z^{(2)}/\nu)\,|\,\mathcal{F}_1]&=\nu^{-2}(C_b+C_W\mathcal{O}_{n_1}^{(1)})
\mathbb{E}[f'_{g}(z^{(2)}/\nu)\,|\,\mathcal{F}_1]\nonumber\\
&=
\mathbb{E}[\nu^{-2}(C_b+C_W\mathcal{O}_{n_1}^{(1)})f'_{g}(z^{(2)}/\nu)\,|\,\mathcal{F}_1],\label{eq:21062023primo}
\end{align}
where the latter equality follows by the $\mathcal{F}_1$-measurability of  $\mathcal{O}_{n_1}^{(1)}$.  
Then
\begin{align}
\mathbb{P}(z^{(2)}\leq y)-\mathbb{P}(z\leq y)=\mathbb E[f'_{g}(z^{(2)}/\nu)(1-\nu^{-2}(C_b+C_W\mathcal{O}_{n_1}^{(1)}))],\quad\text{$y\in\mathbb R$.}\nonumber
\end{align}
Setting
\begin{equation}\label{eq:09052023terzo}
\varphi(n_1):=\nu^{-2}C_W\sqrt{\mathbb{V}\mathrm{ar}(\mathcal{O}_{n_1}^{(1)})}=\nu^{-2}C_W\sqrt{\mathbb{V}\mathrm{ar}(\sigma(z_1^{(1)})^2)}/\sqrt{n_1},
\end{equation}
we have
\begin{align}
\frac{\mathbb{P}(z^{(2)}\leq y)-\mathbb{P}(z\leq y)}{\varphi(n_1)}
=\mathbb{E}\left[f'_{y/\nu_1}(z^{(2)}/\nu)V_{n_1}\right],\label{eq:09052023primoprimo}
\end{align}
where
\begin{align}
V_{n_1}:=\frac{1-\nu^{-2}(C_b+C_W\mathcal{O}_{n_1}^{(1)}))}{\varphi(n_1)}&=-\frac{\mathcal{O}_{n_1}^{(1)}-\mathcal{O}^{(1)}}{\sqrt{\mathbb{V}\mathrm{ar}(\mathcal{O}_{n_1}^{(1)})}}
=-\frac{\sum_{j=1}^{n_1}\sigma(z_j^{(1)})^2-n_1\mathbb E\sigma(z_1^{(1)})^2}{\sqrt{n_1}\sqrt{\mathbb{V}\mathrm{ar}(\sigma(z_1^{(1)})^2)}}
\nonumber
\end{align}
(recall that $\nu$ is defined in \eqref{eq:030720231}).
Taking the modulus in \eqref{eq:09052023primoprimo} and then using that $\|f_g'\|_\infty\leq 1$ uniformly in $y\in\mathbb R$ (see Lemma \ref{le:CGQ}$(i)$) and that
$\mathbb{E}V_{n_1}^2=1$, we have
\[
d_K(z^{(2)},z)\leq\varphi(n_1)\mathbb{E}[|V_{n_1}|]\leq\varphi(n_1),
\]
which,  combined with \eqref{eq:09052023terzo},  proves the statement.\\
\noindent{\it Proof\,\,of\,\,Part\,\,$(ii)$.} Consider the Stein equation \eqref{eq:stein1d}
with $g(w):=\bold{1}_B(\nu w)$,  where $B\subseteq\mathbb R^d$ is a Borel set, and $\nu$ is defined at the beginning of the proof of Part $(i)$.  Let $f_g$ be the unique solution of the Stein equation (see Lemma \ref{le:CGQ}$(ii)$). Then
\begin{align*}
\bold{1}_B(z^{(2)})-\mathbb E\bold{1}_{B}(z)=f'_{g}(z^{(2)}/\nu)-(z^{(2)}/\nu)f_{g}(z^{(2)}/\nu).
\end{align*}
Taking the expectation and arguing as in \eqref{eq:21062023primo},  we have
\[
\mathbb{P}(z^{(2)}\in B)-\mathbb P(z\in B)=
\mathbb{E}[f'_{g}(z^{(2)}/\nu)(1-\nu^{-2}(C_b+C_W\mathcal{O}_{n_1}^{(1)}))].
\]
Along similar computations as for \eqref{eq:09052023primoprimo}, we have
\begin{align*}
\frac{\mathbb{P}(z^{(2)}\in B)-\mathbb P(z\in B)}{\varphi(n_1)}
=\mathbb{E}\left[f'_{g}(z^{(2)}/\nu)V_{n_1}\right].
\end{align*}
Taking the modulus on this relation and then using that $\|f_g'\|_\infty\leq 2$ (see Lemma \ref{le:CGQ}$(ii)$) and that
$\mathbb{E}V_{n_1}^2=1$,  we have
\[
d_{TV}(z^{(2)},z)\leq 2\varphi(n_1)\mathbb{E}[|V_{n_1}|]\leq 2\varphi(n_1),
\]
which,  combined with \eqref{eq:09052023terzo},  proves the statement.\\
\noindent{\it Proof\,\,of\,\,Part\,\,$(iii)$.} Consider the Stein equation \eqref{eq:stein1d}
with $g(y):=h(\nu y)$,  where $h:\mathbb R\to\mathbb R$ is Lipschitz continuous with $\mathrm{Lip}(h)\leq 1$,  and $\nu_1$ is defined at the beginning of the proof of Part $(i)$.  Let $f_g$ be the unique solution of the Stein equation (see Lemma \ref{le:CGQ}$(iii)$). Then
\begin{align*}
h(z^{(2)})-\mathbb E h(z)=f'_{g}(z^{(2)}/\nu)-(z^{(2)}/\nu)f_{g}(z^{(2)}/\nu).
\end{align*}
Taking the expectation and arguing as in \eqref{eq:21062023primo},  we have
\[
\mathbb{E}h(z^{(2)})-\mathbb E h(z)=
\mathbb{E}[f'_{g}(z^{(2)}/\nu)(1-\nu^{-2}(C_b+C_W\mathcal{O}_{n_1}^{(1)}))].
\]
Along similar computations as for \eqref{eq:09052023primoprimo}, we have
\begin{align*}
\frac{\mathbb{E}h(z^{(2)})-\mathbb E h(z)}{\varphi(n_1)}
=\mathbb{E}\left[f'_{g}(z^{(2)}/\nu)V_{n_1}\right].
\end{align*}
Taking the modulus on this relation and then using that $\|f_g'\|_\infty\leq\nu\sqrt{2/\pi}$ (see Lemma \ref{le:CGQ}$(iii)$) and that
$\mathbb{E}V_{n_1}^2=1$,  we have
\[
d_{W_1}(z^{(2)},z)\leq\nu\sqrt{2/\pi}\varphi(n_1)\mathbb{E}[|V_{n_1}|]\leq \nu\sqrt{2/\pi}\varphi(n_1),
\]
which,  combined with \eqref{eq:09052023terzo},  proves the statement.
\\
\noindent$\square$

Note that both Theorem \ref{thm:BFF23} and Theorem \ref{prop:shallow} provide bounds
on $d_s(z^{(2)},z)$,  $s=TV,K,W_1$,  with a common rate $1/\sqrt{n_1}$,  but different constants. 
In Table \ref{tab:shallow-fav-gian} we compare those constants in a special case. 
We observe that the constants given by Theorem \ref{prop:shallow}  are more effective than those in~\cite{BFF23}. 
We also note that Condition \eqref{eq:03072023} is satisfied by the ReLu activation function, while the assumptions of Theorem \ref{thm:BFF23} do not hold for the ReLu.

\begin{table}[ht]
\centering
\small
\begin{tabular}{r|ccc}
 \hline
 & $d_{TV}\left(z^{\left(2\right)},z\right)$ &  $d_{K}\left(z^{\left(2\right)},z\right)$&  $d_{W_1}\left(z^{\left(2\right)},z\right)$ \\ 
  \hline
Thm. \ref{thm:BFF23} & 5.05 & 2.52 & 2.01  \\ 
Thm. \ref{prop:shallow} & 1.68 & 0.84 & 0.67  \\ 
   \hline
\end{tabular}
\caption{Values of the constants given in Theorems  \ref{thm:BFF23} and \ref{prop:shallow}.  Here,  $\sigma(x) = x^3$, 
$\gamma=3$,  $r_2=1$,  $r_1=6$,
$L=1$, $C_b=C_W=1$ and $x=1$. }
\label{tab:shallow-fav-gian} 
\end{table}

\section{A key estimate for the collective observables}\label{sec:collective}

The next theorem provides an estimate for the $L^2$-norm of the random variable
$\mathcal{O}_{\bold{n}_\ell}^{(\ell)}-\mathcal{O}^{(\ell)}$.  In Section~\ref{sec:BOUNDS},
such estimate will play a crucial role in the proofs of the results on the Normal approximation of the output of a deep random Gaussian NN,  both in the convex and in the $1$-Wasserstein distances (see Theorems \ref{thm:convexdistance} and \ref{thm:1Wasserstein}).

Hereon, 
we denote by $\|Y\|_{L^2}:=(\mathbb E[|Y|^2])^{1/2}$
the $L^2$-norm of a real-valued random variable $Y$.

\begin{Theorem}\label{le:fund_est}
Let $\mathrm{GNN}(L,n_0,n_{L+1},\bold{n}_L,\sigma,\bold x,\bold b,\bold W)$ be a deep random Gaussian NN.  Suppose that
the activation function $x\mapsto\sigma(x)$ is such that:\\
\noindent $(i)$ For any $a_1,a_2\geq 0$ and $C_b,C_W>0$,  there exists a polynomial
\[
P(x):=\sum_{k=0}^{m}d_k x^k,
\]
with non-negative coefficients $d_k=d_k(\sigma(\cdot),C_b,C_W)\geq 0$ dependent only on $\sigma(\cdot),C_b, C_W$ and degree
$m\geq 0$ independent of $\sigma(\cdot),a_1,a_2,C_b,C_W$, 
such that
\begin{equation}\label{eq:primasigma}
|\sigma(x\sqrt{C_b+C_W a_2})^2-\sigma(x\sqrt{C_b+C_W a_1})^2|\leq P(|x|)|a_2-a_1|,\quad\text{for all $x\in\mathbb R$.}
\end{equation}
\noindent$(ii)$
For any $\kappa\in\mathbb R$,  $\mathbb E\sigma(\kappa Z)^4<\infty$.
\\

Then, for any $\ell=1,\ldots,L$,  we have
\begin{align*}
\|\mathcal{O}_{\bold{n}_\ell,\sigma^2}^{(\ell)}-\mathcal{O}_{\sigma^2}^{(\ell)}\|_2&\leq\sum_{k=1}^{\ell}(4\sqrt{2}\|P(|Z|)\|_2)^{\ell-k}\frac{c_k}{\sqrt{n_k}}, 
\end{align*}
where
\begin{equation}\label{eq:c20062023}
c_\ell=c_\ell(n_0,\sigma,\bold x,C_b,C_W):=\sqrt{2\mathbb{E}\left[\sigma\left(Z\sqrt{C_b+C_W\mathcal{O}^{(\ell-1)}}\right)^4\right]}<\infty,\quad\text{$\ell=1,\ldots,L$.}
\end{equation}
\end{Theorem}

The proof of the theorem is given later on in this section.  We proceed stating a proposition and a remark,  which clarify the generality of our assumptions on the activation function $\sigma$.  

\begin{Proposition}\label{prop:04052023primo}
The following statements hold:\\
\noindent $(i)$ If $\sigma$ is the perceptron function,  i.e., $\sigma(x):=\bold{1}\{x\geq 0\}$,  $x\in\mathbb R$,  then it satisfies Conditions $(i)$ and $(ii)$ of Theorem \ref{le:fund_est}. \\ 
\noindent $(ii)$ If $\sigma$ is Lipschitz continuous,  then it satisfies Conditions $(i)$ and $(ii)$ of Theorem \ref{le:fund_est}.  In particular,  Condition $(i)$ holds with 
\begin{equation}\label{eq:04052023secondo}
P(x):=\frac{2|\sigma(0)|C_W\mathrm{Lip}(\sigma)}{2\sqrt{C_b}}x+C_W\mathrm{Lip}(\sigma)^2 x^2.
\end{equation}
\noindent $(iii)$ If $\sigma$ is such that $\sigma(\cdot)^2$ is Lipschitz continuous then it satisfies 
Conditions $(i)$ and $(ii)$ of Theorem \ref{le:fund_est}.  In particular,  Condition $(i)$ holds with
\begin{equation}\label{eq:04052023sesto}
P(x)=\frac{\mathrm{Lip}(\sigma^2)C_W}{2\sqrt{C_b}}x.
\end{equation}
\noindent $(iv)$ If $\sigma\in C^1(\mathbb R)$ and $\max\{|\sigma(x)|,|\sigma'(x)|\}\leq r_1+r_2|x|^\gamma$, $x\in\mathbb R$,  for some $r_1,r_2,\gamma\geq 0$, then
$\sigma$ satisfies Conditions $(i)$ and $(ii)$ of Theorem \ref{le:fund_est}.  In particular,  Condition $(i)$ holds with
\[
P(x)=\frac{C_W}{\sqrt{C_b}}x(r_1+r_2|x|^\gamma)
\]
\end{Proposition} 

The proof of the proposition is given later on in this section. 

\begin{Remark}\label{re:sigma}
As a consequence of Proposition \ref{prop:04052023primo},  we have that the most common activation functions satisfy the assumptions of Theorem \ref{le:fund_est}.
Indeed, one can easily prove that the ReLu function $\sigma(x):=x\bold{1}\{x\geq 0\}$,  the sigmoid function $\sigma(x):=(1+\mathrm{e}^{-x})^{-1}$,  the hyperbolic tangent function $\sigma(x):=(\mathrm{e}^{2x}-1)/(\mathrm{e}^{2x}+1)$, the $\sin$ function $\sigma(x):=\sin(x)$,
the softplus function $\sigma(x):=\log(1+\mathrm{e}^{x})$ and the $\mathrm{SWISH}$ function
$\sigma(x):=x/(1+\mathrm{e}^{-x})$ 
are Lipschitz continuous.
We emphasize that the conditions on the activation function of Theorem \ref{le:fund_est} are more general than the one required in \cite{BT}, 
where $\sigma(\cdot)$ is assumed Lipschitz continuous (see Theorem \ref{thm:trevisan} and Proposition \ref{prop:04052023primo}$(ii)$).  We also emphasize that the conditions on the activation function of Theorem \ref{le:fund_est} are
satisfied by the perceptron function which is non-continuous and therefore non-Lipschitz (see Proposition \ref{prop:04052023primo}$(i)$).
Another non-Lipschitz function which satisfies the conditions of Theorem \ref{le:fund_est}
is e.g.  $\sigma(x):=\sqrt{x}\bold{1}\{x\geq 0\}$.  Indeed,  $\sigma^2$ is the ReLu function and therefore Lipschitz continuous (see Proposition \ref{prop:04052023primo}$(iii)$).
\end{Remark}

\noindent{\it Proof\,\,of\,\,Theorem\,\,\ref{le:fund_est}.} We consider separately
the cases of the first hidden layer and that of the following ones.\\
\noindent{\it Case\,\,$\ell=1$.}\\  
Since the random variables $z_i^{(1)}$, $i=1,\ldots,n_1$, are independent and identically distributed with law $\mathcal{N}(0,C_b+C_W\mathcal{O}^{(0)})$,  we have
\begin{align}
\|\mathcal{O}_{n_1}^{(1)}-\mathcal{O}^{(1)}\|_{L^2}
&=
\sqrt{\mathbb{E}\left(\frac{1}{n_{1}}\sum_{j=1}^{n_1}\sigma(z_j^{(1)})^2-
\mathbb{E}\left[\sigma\left(Z\sqrt{C_b+C_W\mathcal{O}^{(0)}}\right)^2\right]\right)^2}\nonumber\\
&=
\sqrt{\mathbb{E}\left(\frac{1}{n_{1}}\sum_{j=1}^{n_1}(\sigma(z_j^{(1)})^2-\mathbb E\sigma(z_1^{(1)})^2)\right)^2}\nonumber\\
&=\frac{1}{\sqrt n_1}\sqrt{\mathbb{V}\mathrm{ar}(\sigma(z_1^{(1)})^2)}\nonumber\\
&\leq\frac{1}{\sqrt n_1}\sqrt{\mathbb{E}\left[\sigma\left(Z\sqrt{C_b+C_W\mathcal{O}^{(0)}}\right)^4\right]}\label{eq:05052023primo}\\
&\leq\frac{c_1}{\sqrt n_1}.\label{eq:25042023dieci}
\end{align}
Note that this latter term is finite due to the assumption $(ii)$. \\
%
\noindent{\it Case\,\,$\ell=2,\ldots,L$.}\\
Take $\ell\in\{2,\ldots,L\}$.
We have already noticed that,   given $\mathcal{F}_{\ell-1}$, 
the random variables $\{z_{i}^{(\ell)}\}_{i=1,\ldots,n_{\ell}}$ are independent with Gaussian law with mean zero and variance 
\[
C_b+C_W\mathcal{O}_{\bold{n}_{\ell-1}}^{(\ell-1)}.
\]
Therefore,  letting $Z_{\ell-1}$ denote a standard Gaussian random variable,  independent of $\mathcal{F}_{\ell-1}$, 
we have
\[
z_{i}^{(\ell)}\overset{d}{=}Z_{\ell-1}\sqrt{C_b+C_W\mathcal{O}_{\bold{n}_{\ell-1}}^{(\ell-1)}},
\quad\text{$i=1,\ldots,n_\ell$}
\]
where the symbold $\overset{d}{=}$ denotes the equality in law (this relation immediately follows computing e.g. the characteristic function of both random variables).
Therefore,
letting $p(\cdot)$ denote the standard Gaussian density, we have 
\begin{equation}\label{eq:25042023primo}
\mathbb{E}[\sigma(z_i^{(\ell)})^r\,|\,\mathcal{F}_{\ell-1}]\overset{d}{=}\int_{\mathbb R}\sigma\left(z\sqrt{C_b+C_W\mathcal{O}_{\bold{n}_{\ell-1}}^{(\ell-1)}}\right)^r p(z)\mathrm{d}z,\quad 
\text{$r\in\{2,4\}$, $i=1,\ldots,n_\ell$}
\end{equation}
and
\begin{align}
\mathbb E\mathcal{O}_{\bold{n}_{\ell}}^{(\ell)}
&=\mathbb{E}\left[\sigma\left(Z_{\ell-1}\sqrt{C_b+C_W\mathcal{O}_{\bold{n}_{\ell-1}}^{(\ell-1)}}\right)^2\right].\nonumber
\end{align}
By assumption $(i)$ and the fact that $Z_{\ell-1}$ is independent of $\mathcal{F}_{\ell-1}$, we have
\begin{align}
|\mathbb E\mathcal{O}_{\bold{n}_{\ell}}^{(\ell)}-\mathcal{O}^{(\ell)}|&=\Big|\mathbb{E}\left[\sigma\left(Z_{\ell-1}\sqrt{C_b+C_W\mathcal{O}_{\bold{n}_{\ell-1}}^{(\ell-1)}}\right)^2\right]-\mathbb{E}\left[\sigma\left(Z_{\ell-1}\sqrt{C_b+C_W\mathcal{O}^{(\ell-1)}}\right)^2\right]\Big|\nonumber\\
&\leq\mathbb{E}P(|Z_{\ell-1}|)|\mathcal{O}_{\bold{n}_{\ell-1}}^{(\ell-1)}-\mathcal{O}^{(\ell-1)}|
\nonumber\\
&\leq\mathbb{E}P(|Z|)\|\mathcal{O}_{\bold{n}_{\ell-1}}^{(\ell-1)}-\mathcal{O}^{(\ell-1)}\|_{L^2}.\nonumber
\end{align}
Therefore
\begin{align}
\|\mathcal{O}_{\bold{n}_{\ell}}^{(\ell)}-\mathcal{O}^{(\ell)}\|_{L^2}&\leq\|\mathcal{O}_{\bold{n}_{\ell}}^{(\ell)}-\mathbb{E}\mathcal{O}_{\bold{n}_{\ell}}^{(\ell)}\|_{L^2}+|\mathbb{E}\mathcal{O}_{\bold{n}_{\ell}}^{(\ell)}-\mathcal{O}^{(\ell)}|\nonumber\\
&=\sqrt{\mathbb{V}\mathrm{ar}(\mathcal{O}_{\bold{n}_{\ell}}^{(\ell)})}
+|\mathbb{E}\mathcal{O}_{\bold{n}_{\ell}}^{(\ell)}-\mathcal{O}^{(\ell)}|
\nonumber\\
&\leq\sqrt{\mathbb{V}\mathrm{ar}(\mathcal{O}_{\bold{n}_{\ell}}^{(\ell)})}+
\mathbb E P(|Z|)\|\mathcal{O}_{\bold{n}_{\ell-1}}^{(\ell-1)}-\mathcal{O}^{(\ell-1)}\|_{L^2}.
\label{eq:24042023primo}
\end{align}
Note that
\begin{align}
&\mathbb{V}\mathrm{ar}(\mathcal{O}_{\bold{n}_{\ell}}^{(\ell)})=
\frac{1}{n_{\ell}^2}\left(\mathbb E\left(
\sum_{j=1}^{n_{\ell}}\sigma(z_j^{(\ell)})^2\right)^2
-n_{\ell}^2(\mathbb E \sigma(z_1^{(\ell)})^2)^2
\right)\nonumber\\
&=
\frac{1}{n_{\ell}^2}\left(n_{\ell}\mathbb E \sigma(z_1^{(\ell)})^4+n_{\ell}(n_{\ell}-1)\mathbb{E}\sigma(z_1^{(\ell)})^2\sigma(z_2^{(\ell)})^2
-n_{\ell}^2(\mathbb E \sigma(z_1^{(\ell)})^2)^2
\right)\nonumber\\
&=
\frac{1}{n_{\ell}^2}\left(n_{\ell}\mathbb E \sigma(z_1^{(\ell)})^4-n_{\ell}\mathbb{E}\sigma(z_1^{(\ell)})^2\sigma(z_2^{(\ell)})^2
+n_{\ell}^2\mathbb{C}\mathrm{ov}(\sigma(z_1^{(\ell)})^2,\sigma(z_2^{(\ell)})^2)
\right)\nonumber\\
&=
\frac{1}{n_{\ell}^2}\left(n_{\ell}\mathbb{V}\mathrm{ar}(\sigma(z_1^{(\ell)})^2)
+(n_{\ell}^2-n_{\ell})\mathbb{C}\mathrm{ov}(\sigma(z_1^{(\ell)})^2,\sigma(z_2^{(\ell)})^2)\right)\nonumber\\
&=
\frac{1}{n_{\ell}}\mathbb{V}\mathrm{ar}(\sigma(z_1^{(\ell)})^2)
+\left(1-\frac{1}{n_{\ell}}\right)\mathbb{C}\mathrm{ov}(\sigma(z_1^{(\ell)})^2,\sigma(z_2^{(\ell)})^2).\label{eq:25042023seconda}
\end{align}
By \eqref{eq:25042023primo} we have
\[
\mathbb{E}[\sigma(z_1^{(\ell)})^4]=\mathbb E\sigma\left(Z_{\ell-1}\sqrt{C_b+C_W\mathcal{O}_{\bold{n}_{\ell-1}}^{(\ell-1)}}\right)^4.
\]
By assumption $(i)$,  we have
\begin{equation*}
\sigma\left(Z_{\ell-1}\sqrt{C_b+C_W \mathcal{O}_{\bold{n}_{\ell-1}}^{(\ell-1)}}\right)^2\leq P(|Z_{\ell-1}|)|\mathcal{O}_{\bold{n}_{\ell-1}}^{(\ell-1)}-\mathcal{O}^{(\ell-1)}|
+\sigma\left(Z_{\ell-1}\sqrt{C_b+C_W \mathcal{O}^{(\ell-1)}}\right)^2,\quad\text{$\mathbb P$-a.s.}
\end{equation*}
and so (using that $Z_{\ell-1}$ is independent of $\mathcal{F}_{\ell-1}$ and the inequality $(a+b)^2\leq 2a^2+2b^2$,  $a,b\in\mathbb R$)
\begin{align}
\mathbb{V}\mathrm{ar}(\sigma(z_1^{(\ell)})^2)&\leq\mathbb{E}[\sigma(z_1^{(\ell)})^4]\nonumber\\
&\leq 2\mathbb E P(|Z|)^{2}\|\mathcal{O}_{\bold{n}_{\ell-1}}^{(\ell-1)}-\mathcal{O}^{(\ell-1)}\|_2^2
+2\mathbb E\sigma\left(Z\sqrt{C_b+C_W\mathcal{O}^{(\ell-1)}}\right)^4
\nonumber\\
&=A\|\mathcal{O}_{\bold{n}_{\ell-1}}^{(\ell-1)}-\mathcal{O}^{(\ell-1)}\|_{L^2}^2
+c_\ell^2,
\label{eq:25042023terzo}
\end{align}
where $A:=2\mathbb E P(|Z|)^2<\infty$.  Note that the quantities $\mathcal{O}^{(\ell)}$,  $\ell=2,\ldots,L$,
are all finite due to the assumption $(ii)$. 
Then,  again by assumption $(ii)$,  we have
\[
\mathbb E\sigma\left(Z\sqrt{C_b+C_W\mathcal{O}^{(\ell-1)}}\right)^4<\infty,\quad\text{for any $\ell=2,\ldots,L$}
\]
and therefore $c_\ell<\infty$, for any $\ell=2,\ldots,L$.
By the conditional independence of the random variables
$z_1^{(\ell)}$ and $z_2^{(\ell)}$, given $\mathcal{F}_{\ell-1}$, we have
\begin{align}
&\mathbb{C}\mathrm{ov}(\sigma(z_1^{(\ell)})^2,\sigma(z_2^{(\ell)})^2)=
\mathbb{C}\mathrm{ov}(\mathbb{E}[\sigma(z_1^{(\ell)})^2\,|\,\mathcal{F}_{\ell-1}],\mathbb{E}[\sigma(z_2^{(\ell)})^2\,|\,\mathcal{F}_{\ell-1}])\nonumber\\
&\qquad\qquad
=\mathbb E[(\mathbb{E}[\sigma(z_1^{(\ell)})^2\,|\,\mathcal{F}_{\ell-1}]-\mathbb  E\sigma(z_1^{(\ell)})^2)(\mathbb{E}[\sigma(z_2^{(\ell)})^2\,|\,\mathcal{F}_{\ell-1}]-\mathbb E\sigma(z_2^{(\ell)})^2)],\nonumber
\end{align}
and so by the Cauchy-Schwarz inequality and \eqref{eq:25042023primo} we have
\begin{equation}\label{eq:15092023primo}
|\mathbb{C}\mathrm{ov}(\sigma(z_1^{(\ell)})^2,\sigma(z_2^{(\ell)})^2)|\leq
\mathbb E[(\mathbb{E}[\sigma(z_1^{(\ell)})^2\,|\,\mathcal{F}_{\ell-1}]-\mathbb  E\sigma(z_1^{(\ell)})^2)^2].
\end{equation}
Letting $\mathbb P_X$ denote the law of a random variable $X$ and
using again \eqref{eq:25042023primo}, we have that the random variable 
$\mathbb{E}[\sigma(z_1^{(\ell)})^2\,|\,\mathcal{F}_{\ell-1}]-\mathbb  E\sigma(z_1^{(\ell)})^2$ has the same law as the random variable
\begin{align}
&\int_{\mathbb R}\left(\sigma\left(z\sqrt{C_b+C_W\mathcal{O}_{\bold{n}_{\ell-1}}^{(\ell-1)}}\right)^2-\int_{[0,\infty)}\sigma\left(z\sqrt{C_b+C_Wy}\right)^2\mathbb{P}_{\mathcal{O}_{\bold{n}_{\ell-1}}^{(\ell-1)}}(\mathrm{d}y)\right)p(z)\mathrm{d}z\nonumber\\
&\quad=
\int_{[0,\infty)\times\mathbb R}\left(\sigma\left(z\sqrt{C_b+C_W\mathcal{O}_{\bold{n}_{\ell-1}}^{(\ell-1)}}\right)^2-\sigma\left(z\sqrt{C_b+C_Wy}\right)^2\right)\mathbb{P}_{\mathcal{O}_{\bold{n}_{\ell-1}}^{(\ell-1)}}(\mathrm{d}y)p(z)\mathrm{d}z.\nonumber
\end{align}
Therefore,  by \eqref{eq:15092023primo} and Jensen's inequality we have
\begin{align}
&|\mathbb{C}\mathrm{ov}(\sigma(z_1^{(\ell)})^2,\sigma(z_2^{(\ell)})^2)|
\nonumber\\
&\qquad
\leq\mathbb{E}
\int_{[0,\infty)\times\mathbb R}\left(\sigma\left(z\sqrt{C_b+C_W\mathcal{O}_{\bold{n}_{\ell-1}}^{(\ell-1)}}\right)^2-\sigma\left(z\sqrt{C_b+C_Wy}\right)^2\right)^2 \mathbb{P}_{\mathcal{O}_{\bold{n}_{\ell-1}}^{(\ell-1)}}(\mathrm{d}y)p(z)\mathrm{d}z
\nonumber
\end{align}
By assumption $(i)$ it then follows that
%
%
%
\begin{align}
|\mathbb{C}\mathrm{ov}(\sigma(z_1^{(\ell)})^2,\sigma(z_2^{(\ell)})^2)|
&\leq(\mathbb E P(|Z|))^2\int_{[0,\infty)}\mathbb E|\mathcal{O}_{\bold{n}_{\ell-1}}^{(\ell-1)}-y|^2\mathbb{P}_{\mathcal{O}_{\bold{n}_{\ell-1}}^{(\ell-1)}}(\mathrm{d}y)
\leq A\mathbb{V}\mathrm{ar}\left(\mathcal{O}_{\bold{n}_{\ell-1}}^{(\ell-1)}\right),\label{eq:25042023quarto}
\end{align}
where the latter relation 
follows by 
the definition of the constant $A$.
Combining \eqref{eq:25042023seconda},  \eqref{eq:25042023terzo} and \eqref{eq:25042023quarto},  we have
\begin{align}
\mathbb{V}\mathrm{ar}(\mathcal{O}_{\bold{n}_\ell}^{(\ell)})
&\leq A(\mathbb{V}\mathrm{ar}(\mathcal{O}_{\bold{n}_{\ell-1}}^{(\ell-1)})+\|\mathcal{O}_{\bold{n}_{\ell-1}}^{(\ell-1)}-\mathcal{O}_{\sigma^2}^{(\ell-1)}\|_{L^2}^2)
+\frac{c_{\ell}^2}{n_\ell}.\nonumber
\end{align}
Iterating this inequality we have
\begin{align}
\mathbb{V}\mathrm{ar}(\mathcal{O}_{\bold{n}_\ell}^{(\ell)})
&\leq A\left(A(\mathbb{V}\mathrm{ar}(\mathcal{O}_{\bold{n}_{\ell-2}}^{(\ell-2)})+
\|\mathcal{O}_{\bold{n}_{\ell-2}}^{(\ell-2)}-\mathcal{O}^{(\ell-2)}\|_{L^2}^2)
+\frac{c_{\ell-1}^2}{n_{\ell-1}}\right)\nonumber\\
&\qquad\qquad\qquad\qquad\qquad
+A\|\mathcal{O}_{\bold{n}_{\ell-1}}^{(\ell-1)}-\mathcal{O}^{(\ell-1)}\|_{L^2}^2
+\frac{c_{\ell}^2}{n_\ell}\nonumber\\
&=A^2\mathbb{V}\mathrm{ar}(\mathcal{O}_{\bold{n}_{\ell-2}}^{(\ell-2)})
+A^2
\|\mathcal{O}_{\bold{n}_{\ell-2}}^{(\ell-2)}-\mathcal{O}^{(\ell-2)}\|_{L^2}^2
+A\|\mathcal{O}_{\bold{n}_{\ell-1}}^{(\ell-1)}-\mathcal{O}^{(\ell-1)}\|_{L^2}^2\nonumber\\
&\qquad\qquad\qquad\qquad\qquad
+A\frac{c_{\ell-1}^2}{n_{\ell-1}}+\frac{c_{\ell}^2}{n_\ell}\nonumber\\
&\leq
A^3\mathbb{V}\mathrm{ar}(\mathcal{O}_{\bold{n}_{\ell-3}}^{(\ell-3)})
+A^3
\|\mathcal{O}_{\bold{n}_{\ell-3}}^{(\ell-3)}-\mathcal{O}^{(\ell-3)}\|_{L^2}^2
+A^2
\|\mathcal{O}_{\bold{n}_{\ell-2}}^{(\ell-2)}-\mathcal{O}^{(\ell-2)}\|_{L^2}^2\nonumber\\
&\qquad\qquad\qquad
+A\|\mathcal{O}_{\bold{n}_{\ell-1}}^{(\ell-1)}-\mathcal{O}^{(\ell-1)}\|_{L^2}^2
+A^2\frac{c_{\ell-2}^2}{n_{\ell-2}}+A\frac{c_{\ell-1}^2}{n_{\ell-1}}+\frac{c_{\ell}^2}{n_\ell}\nonumber\\
&\leq
A^{\ell-1}\mathbb{V}\mathrm{ar}(\mathcal{O}_{n_1}^{(1)})
+\sum_{k=1}^{\ell-1}A^{\ell-k}
\|\mathcal{O}_{\bold{n}_{k}}^{(k)}-\mathcal{O}^{(k)}\|_{L^2}^2
+\sum_{k=2}^{\ell}A^{\ell-k}\frac{c_{k}^2}{n_{k}}\nonumber\\
&=
2A^{\ell-1}\mathbb{V}\mathrm{ar}(\mathcal{O}_{n_1}^{(1)})
+\sum_{k=2}^{\ell-1}A^{\ell-k}
\|\mathcal{O}_{\bold{n}_{k}}^{(k)}-\mathcal{O}^{(k)}\|_{L^2}^2
+\sum_{k=2}^{\ell}A^{\ell-k}\frac{c_{k}^2}{n_{k}},\nonumber
\end{align}
for any $\ell=2,\ldots,L$,
where the latter equality follows noticing that $\mathbb E \mathcal{O}_{n_1}^{(1)}=\mathcal{O}^{(1)}$.
Here, we adopt the usual convention $\sum_{k=k_1}^{k_2}\cdots=0$ if $k_1>k_2$.
Note that by \eqref{eq:05052023primo} we have
\[
\mathbb{V}\mathrm{ar}(\mathcal{O}_{n_1}^{(1)})\leq
\frac{1}{n_1}\mathbb{E}\sigma\left(Z\sqrt{C_b+C_W\mathcal{O}^{(0)}}\right)^4,
\]
and so
\[
2\mathbb{V}\mathrm{ar}(\mathcal{O}_{n_1}^{(1)})\leq\frac{c_1^2}{n_1}.
\]
Consequently,
\[
\mathbb{V}\mathrm{ar}(\mathcal{O}_{\bold{n}_\ell}^{(\ell)})
\leq\sum_{k=2}^{\ell-1}A^{\ell-k}
\|\mathcal{O}_{\bold{n}_{k}}^{(k)}-\mathcal{O}^{(k)}\|_{L^2}^2
+\sum_{k=1}^{\ell}A^{\ell-k}\frac{c_{k}^2}{n_{k}}.\nonumber
\]
Combining this inequality with \eqref{eq:24042023primo} and using the elementary relation
$\sqrt{a_1+a_2}\leq\sqrt a_1+\sqrt a_2$,  $a_1,a_2\geq 0$,
for any $\ell=2,\ldots,L$, we have
\begin{align}
\|\mathcal{O}_{\bold{n}_{\ell}}^{(\ell)}-\mathcal{O}^{(\ell)}\|_{L^2}&\leq
\sum_{k=2}^{\ell-1}A^{(\ell-k)/2}
\|\mathcal{O}_{\bold{n}_{k}}^{(k)}-\mathcal{O}^{(k)}\|_{L^2}
+\sum_{k=1}^{\ell}A^{(\ell-k)/2}\frac{c_{k}}{\sqrt{n_{k}}}\nonumber\\
&\qquad\qquad\qquad
+\mathbb E P(|Z|)\|\mathcal{O}_{\bold{n}_{\ell-1}}^{(\ell-1)}-\mathcal{O}^{(\ell-1)}\|_{L^2}\nonumber\\
&\leq
\sum_{k=2}^{\ell-1}B^{(\ell-k)/2}
\|\mathcal{O}_{\bold{n}_{k}}^{(k)}-\mathcal{O}^{(k)}\|_{L^2}
+\sum_{k=1}^{\ell}B^{(\ell-k)/2}\frac{c_{k}}{\sqrt{n_{k}}},\label{eq:05052023secondo}
\end{align}
where we used that $\mathbb E P(|Z|)\leq A^{1/2}$, and we set $B:=4A$.  Recalling that $\sqrt{A}=\sqrt{2}\|P(|Z|)\|_2$ and the definition of $B$, one easily realizes that the claim reads as 
\begin{equation}\label{eq:16092023primo}
\|\mathcal{O}_{\bold{n}_{2}}^{(\ell)}-\mathcal{O}^{(\ell)}\|_{L^2}\leq
\sum_{j=1}^{\ell}(2\sqrt{B})^{\ell-j}\frac{c_j}{\sqrt{n_j}},\quad\text{$\ell=2,\ldots,L$. }
\end{equation} 
Now we use the relation \eqref{eq:05052023secondo} to prove \eqref{eq:16092023primo} by induction on $\ell=2,\ldots,L$.  Taking $\ell=2$ in \eqref{eq:05052023secondo} we have
\[
\|\mathcal{O}_{\bold{n}_{2}}^{(2)}-\mathcal{O}^{(2)}\|_{L^2}
\leq B^{1/2}\frac{c_1}{\sqrt{n_1}}+\frac{c_2}{\sqrt{n_2}}\leq 2B^{1/2}\frac{c_1}{\sqrt{n_1}}+\frac{c_2}{\sqrt{n_2}},
\] 
i.e.,  \eqref{eq:16092023primo} with $\ell=2$. Now,  suppose that 
\[
\|\mathcal{O}_{\bold{n}_{k}}^{(k)}-\mathcal{O}^{(k)}(\bold x)\|_{L^2}
\leq\sum_{j=1}^{k}(2\sqrt{B})^{k-j}\frac{c_j}{\sqrt{n_j}},\quad\text{for any $k=2,\ldots,\ell-1$.}
\] 
Then by \eqref{eq:05052023secondo},  we have
\begin{align}
\|\mathcal{O}_{\bold{n}_{\ell}}^{(\ell)}-\mathcal{O}^{(\ell)}\|_{L^2}&\leq
\sum_{k=2}^{\ell-1}B^{(\ell-k)/2}
\sum_{j=1}^{k}2^{k-j}B^{(k-j)/2}\frac{c_j}{\sqrt{n_j}}
+\sum_{k=1}^{\ell}B^{(\ell-k)/2}\frac{c_{k}}{\sqrt{n_{k}}}\nonumber\\
&\leq
\sum_{k=2}^{\ell-1}\sum_{j=1}^{k}2^{k-j}B^{(\ell-j)/2}\frac{c_j}{\sqrt{n_j}}
+\sum_{k=1}^{\ell}B^{(\ell-k)/2}\frac{c_{k}}{\sqrt{n_{k}}}\nonumber\\
&\leq
\sum_{j=1}^{\ell-1}\sum_{k=j}^{\ell-1}2^{k-j}B^{(\ell-j)/2}\frac{c_j}{\sqrt{n_j}}
+\sum_{k=1}^{\ell}B^{(\ell-k)/2}\frac{c_{k}}{\sqrt{n_{k}}}\nonumber\\
&=
\sum_{j=1}^{\ell-1}\left(\sum_{k=j}^{\ell-1}2^{k-j}+1\right)B^{(\ell-j)/2}\frac{c_j}{\sqrt{n_j}}
+\frac{c_{\ell}}{\sqrt{n_{\ell}}}\nonumber\\
&=
\sum_{j=1}^{\ell-1}(2\sqrt{B})^{\ell-j}\frac{c_j}{\sqrt{n_j}}
+\frac{c_{\ell}}{\sqrt{n_{\ell}}}=\sum_{j=1}^{\ell}(2\sqrt{B})^{\ell-j}\frac{c_j}{\sqrt{n_j}}
,\nonumber
\end{align}
where we used that
\[
\sum_{k=j}^{\ell-1}2^{k-j}+1=\sum_{s=0}^{\ell-j-1}2^s+1=2^{\ell-j},\quad\text{for any $j=1,\ldots,\ell-1$}
\]
since $\{2^s\}_{s\geq 0}$ is a geometric progression. The proof is completed. 
\\
\noindent$\square$

\noindent{\it Proof\,\,of\,\,Proposition\,\,\ref{prop:04052023primo}.}\\
\noindent{\it Proof\,\,of\,\,$(i)$.} The claim immediately follows noticing that,  due to the positivity of the quantities $\sqrt{C_b+C_W a_2}$ and $\sqrt{C_b+C_W a_1}$,  we have
$|\sigma(x\sqrt{C_b+C_W a_2})^2-\sigma(x\sqrt{C_b+C_W a_1})^2|=0$,  for any $x\in\mathbb R$.\\
\noindent{\it Proof\,\,of\,\,$(ii)$.} Since $\sigma$ is Lipschitz continuous,  we have
\begin{align}
&|\sigma(x\sqrt{C_b+C_W a_2})^2-\sigma(x\sqrt{C_b+C_W a_1})^2|\nonumber\\
&=|\sigma(x\sqrt{C_b+C_W a_2})-\sigma(x\sqrt{C_b+C_W a_1})|
|\sigma(x\sqrt{C_b+C_W a_2})+\sigma(x\sqrt{C_b+C_W a_1})|\nonumber\\
&\leq\mathrm{Lip}(\sigma)|x||\sqrt{C_b+C_W a_2}-\sqrt{C_b+C_W a_1}|(|\sigma(x\sqrt{C_b+C_W a_2})|+|\sigma(x\sqrt{C_b+C_W a_1})|)\nonumber\\
&\leq\mathrm{Lip}(\sigma)|x||\sqrt{C_b+C_W a_2}-\sqrt{C_b+C_W a_1}|[2|\sigma(0)|+\mathrm{Lip}(\sigma)|x|(\sqrt{C_b+C_W a_2}+\sqrt{C_b+C_W a_1})],
\label{eq:04052023terzo}
\end{align}
where the latter inequality follows noticing that (again by the Lipschitz continuity of $\sigma$),  for any $x,\kappa\in\mathbb R$,
\begin{equation}\label{eq:04052023quinto}
|\sigma(\kappa x)|\leq |\sigma(0)|+\mathrm{Lip}(\sigma(\cdot))|\kappa||x|.
\end{equation}
Multiplying and dividing the term in \eqref{eq:04052023terzo} by  $\sqrt{C_b+C_W a_2}+\sqrt{C_b+C_W a_1}$,  we easily have that
\begin{align}
&|\sigma(x\sqrt{C_b+C_W a_2})^2-\sigma(x\sqrt{C_b+C_W a_1})^2|\nonumber\\
&\leq C_W\mathrm{Lip}(\sigma)|x|\left(\frac{2|\sigma(0)|}{\sqrt{C_b+C_W a_2}+\sqrt{C_b+C_W a_1}}+\mathrm{Lip}(\sigma)|x|\right)|a_2-a_1|\nonumber\\
&\leq P(|x|)|a_2-a_1|,\nonumber
\end{align}
where $P(x)$ is defined by \eqref{eq:04052023secondo}.  As far as Condition $(ii)$ of Theorem \ref{le:fund_est} is concerned,  we note that by \eqref{eq:04052023quinto}
\[
\sigma(\kappa Z)^4\leq(|\sigma(0)|+\mathrm{Lip}(\sigma)|\kappa||Z|)^4,
\quad\text{almost surely}
\]
and so the claim is an immediate consequence of the fact that the absolute moments of $Z$ are finite.\\
\noindent{\it Proof\,\,of\,\,$(iii)$.}
If $\sigma(\cdot)^2$ is Lipschitz continuous,  then
\begin{align}
&|\sigma(x\sqrt{C_b+C_W a_2})^2-\sigma(x\sqrt{C_b+C_W a_1})^2|\nonumber\\
&\quad\quad
\leq\mathrm{Lip}(\sigma^2)|x||\sqrt{C_b+C_W a_2}-\sqrt{C_b+C_W a_1}|\nonumber\\
&\quad\quad
=\mathrm{Lip}(\sigma^2)|x||\sqrt{C_b+C_W a_2}-\sqrt{C_b+C_W a_1}|
\frac{\sqrt{C_b+C_W a_2}+\sqrt{C_b+C_W a_1}}{\sqrt{C_b+C_W a_2}+\sqrt{C_b+C_W a_1}}\nonumber\\
&\quad\quad
=\mathrm{Lip}(\sigma^2)|x||(C_b+C_W a_2)-(C_b+C_W a_1)|
\frac{1}{\sqrt{C_b+C_W a_2}+\sqrt{C_b+C_W a_1}}\nonumber\\
&\quad\quad
\leq\frac{\mathrm{Lip}(\sigma^2)C_W}{2\sqrt{C_b}}|x||a_2-a_1|,\nonumber
\end{align}
which shows that Condition $(i)$ of Theorem \ref{le:fund_est} holds with $P(x)$ given by \eqref{eq:04052023sesto}.
Furthermore,  using \eqref{eq:04052023quinto} with $\sigma(\cdot)^2$ in place of $\sigma(\cdot)$,  for any $\kappa\in\mathbb R$,
\[
\sigma(\kappa Z)^4\leq\left(|\sigma(0)|^2+\mathrm{Lip}(\sigma(\cdot)^2)|\kappa||Z|\right)^2,
\quad\text{almost surely.}
\]
Therefore,  Condition $(ii)$ of Theorem \ref{le:fund_est} 
is an immediate consequence of this latter relation and the fact that the absolute moments of $Z$ are finite.\\
\noindent{\it Proof\,\,of\,\,$(iv)$.} Condition $(ii)$ of Theorem \ref{le:fund_est} immediately follows noticing that,  for any $\kappa\in\mathbb R$,
\[
\mathbb E\sigma(\kappa Z)^4\leq\mathbb E(r_1+r_2|\kappa|^\gamma|Z|^\gamma)^4<\infty.
\]
As far as Condition $(i)$ is concerned,  note that by the mean value theorem, for any $x\in\mathbb R$,  there exists $\xi\in\left(\min\{x\sqrt{C_b+C_W a_1},x\sqrt{C_b+C_W a_1}\},\max\{x\sqrt{C_b+C_W a_1},x\sqrt{C_b+C_W a_1}\}\right)$ such that
\begin{align}
|\sigma(x\sqrt{C_b+C_W a_2})^2-\sigma(x\sqrt{C_b+C_W a_1})^2|&=2|\sigma'(\xi)||x|
|\sqrt{C_b+C_W a_2}-\sqrt{C_b+C_W a_1}|\nonumber\\
&=\frac{2C_W|\sigma'(\xi)|}{\sqrt{C_b+C_W a_2}+\sqrt{C_b+C_W a_1}}|x|
|a_2-a_1|.\nonumber
\end{align}
The claim follows noticing that by assumption $|\sigma'(x)|\leq r_1+r_2|x|^\gamma$, for any $x\in\mathbb R$.
\\
\noindent$\square$

\section{Normal approximation of deep random Gaussian NNs}\label{sec:BOUNDS}

\subsection{Normal approximation of deep random Gaussian NNs in the convex distance}

The following theorem holds.

\begin{Theorem}\label{thm:convexdistance}
Let $\mathrm{GNN}(L,n_0,n_{L+1},\bold{n}_{L},\sigma;\bold x;\bold b,\bold W)$ be a deep random Gaussian NN,  and let
the notation and assumptions of Theorem \ref{le:fund_est} prevail.  Then
\begin{align}
d_c(\bold{z}^{(L+1)},\bold{z})&\leq C_1\left(\sum_{k=1}^{L}[4\sqrt{2}\|P(|Z|)\|_{L^2}]^{L-k}\frac{c_k}{\sqrt{n_k}}\right)\nonumber\\
&\leq C_i\left(\sum_{k=1}^{L}[4\sqrt{2}\|P(|Z|)\|_{L^2}]^{L-k}\frac{c_k}{\sqrt{n_k}}\right),\quad\text{$i=2,3$}\nonumber
\end{align}
where the constants $c_k$, $k=1,\ldots,L$,  are defined by \eqref{eq:c20062023},
\begin{equation}\label{eq: C_1 esplicita}
C_1=C_1(n_0,n_{L+1},\bold x, \sigma,C_b,C_W):=C_W\left(\frac{80}{(C_b+C_W\mathcal{O}^{(L+1)})^{3/2}}+\frac{48}{C_b+C_W\mathcal{O}^{(L+1)}}+20\sqrt{2}\right)n_{L+1}^{59/24},
\end{equation}
\[
C_2=C_2(n_{L+1},C_b,C_W):=C_W\left(\frac{80}{C_b^{3/2}}+\frac{48}{C_b}+20\sqrt{2}\right)n_{L+1}^{59/24}
\]
and
\[
C_3=C_3(n_0,n_{L+1},\bold x, \sigma,C_W):=C_W\left(\frac{80}{(C_W\mathcal{O}_{\sigma^2}^{(L+1)})^{3/2}}+\frac{48}{C_W\mathcal{O}_{\sigma^2}^{(L+1)}}+20\sqrt{2}\right)n_{L+1}^{59/24}.
\]
\end{Theorem}

\begin{Remark}\label{re:18092023secondo} Let $\mathrm{GNN}(L,n_0,n_{L+1},\bold{n}_{L},\sigma;\bold x;\bold b,\bold W)$ be a deep random Gaussian NN.  Theorem 3.5 in \cite{FHMNP} shows that if 
\begin{equation}\label{eq:19092023}
c_2 n\leq n_1,\ldots,n_L\leq c_1 n,\quad\text{for some constants $c_1\geq c_2>0$ and some $n\geq 1$}
\end{equation}
and $\sigma$ is polynomially bounded to order $r\geq 1$,  then there exists a constant $C_0$ such that
\[
d_c(\bold{z}^{(L+1)},\bold{z})\leq C_0 n^{-1/2}.
\]
We stress that the bounds in Theorem \ref{thm:convexdistance} provide (under different assumptions on $\sigma$ and without any condition on the widths of the hidden layers) a detailed description of the analytical dependence of the upper estimate on the parameters of the model. 
\end{Remark}

\noindent{\it Proof.} Throughout this proof,  for ease of notation, we put
\[
\kappa:=d_c(\bold{z}^{(L+1)},\bold{z})\quad\text{and}\quad\gamma:=C_W\|\mathcal{O}_{\bold{n}_L}^{(L)}-\mathcal{O}^{(L)}\|_{L^2}.
\]
We preliminary note that it suffices to prove that
\begin{align}
\kappa\leq(80\|\Sigma_{n_{L+1}}^{-1}\|_{op}^{3/2}+48\|\Sigma_{n_{L+1}}^{-1}\|_{op}+20\sqrt{2})n_{L+1}^{59/24}\gamma.
\label{eq:02052023undici}
\end{align} 
Indeed,  the claim then follows by Theorem \ref{le:fund_est},  noticing that
\begin{equation}\label{eq:10052023secondo}
\|\Sigma_{n_{L+1}}^{-1}\|_{op}=\frac{1}{C_b+C_W\mathcal{O}^{(L)}}
\end{equation}
and that $C_1\leq C_i$, $i=2,3$.

If $\gamma>1/\mathrm{e}$,  then the inequality \eqref{eq:02052023undici} holds since $\kappa\leq 1$ (which follows by the definition of the convex distance) and 
\[
20\sqrt{2}n_{L+1}^{59/24}\gamma>20\sqrt{2}\gamma>20\sqrt{2}/3>1.
\]

From now on,  we assume $\gamma\leq1/\mathrm{e}$.
Let $h\in\mathcal{I}_{n_{L+1}}$ (i.e.,  $h$ is the indicator function of a measurable convex set in $\mathbb R^{n_{L+1}}$) be arbitrarily fixed and let 
\[
h_{t}(\bold{y}):=\mathbb E h(\sqrt{t}\bold{z}+\sqrt{1-t}\bold{y}),\quad t\in (0,1),\quad\bold y\in\mathbb R^{n_{L+1}}.
\] 
For any $t\in (0,1)$,  by Lemma \ref{le:Stein2}$(i)$
\[
f_{t,h}(\bold y):=\frac{1}{2}\int_t^1\frac{1}{1-s}(\mathbb E h(\sqrt{t}\bold{z}+\sqrt{1-t}\bold{y})-\mathbb E h(\bold{z}))\mathrm{d}s,
\quad\bold y\in\mathbb{R}^{n_{L+1}}
\]
solves the Stein equation \eqref{eq:SteinNew} with $n_{L+1}$,  $h_{t}$, $\Sigma_{n_{L+1}}$ defined by \eqref{eq:09052023primo}, and $\bold z$,
in place of $d$,  $g$,  $\Sigma$ and $\bold{N}_\Sigma$,  respectively, i.e.,
\begin{equation}\label{eq:09052023secondo}
h_{t}(\bold y)-\mathbb{E}[h_{t}(\bold z)]=\langle\bold y,\nabla f_{t,h}(\bold y)\rangle_{n_{L+1}}-\langle\Sigma_{n_{L+1}},\mathrm{Hess}\,f_{t,h}(\bold y)\rangle_{H.S.},\quad\bold y\in\mathbb R^{n_{L+1}}.
\end{equation}
By Lemma \ref{le:smoothingdconvex} it follows
\begin{equation}\label{eq:02052023secondo}
\kappa\leq\frac{4}{3}\sup_{h\in\mathcal{I}_{n_{L+1}}}|\mathbb E h_{t}(\bold{z}^{(L+1)})-\mathbb E h_{t}(\bold{z})|+\frac{20 n_{L+1}}{\sqrt 2}\frac{\sqrt t}{1-t},\quad t\in (0,1).
\end{equation}
Without loss of generality,  hereafter we assume that $\bold{z}$ is independent of $\mathcal{F}_L$. Therefore,  by \eqref{eq:09052023secondo} we have
\begin{align}
\mathbb E [h_{t}(\bold{z}^{(L+1)})-h_{t}(\bold{z})\,|\,\mathcal{F}_L]
&=\sum_{i=1}^{n_{L+1}}\mathbb E[z_i^{(L+1)}\partial_{i}f_{t,h}(\bold{z}^{(L+1)})\,|\,\mathcal{F}_L]
-\sum_{i,j=1}^{n_{L+1}}\Sigma_{n_{L+1}}(i,j)\mathbb{E}[\partial_{ij}^2 f_{t,h}(\bold{z}^{(L+1)})\,|\,\mathcal{F}_L].
\label{eq:09052023quarto}
\end{align}
By Lemma \ref{le:Stein2}$(i)$,
for any $t\in (0,1)$,  the mapping $\bold y\mapsto\partial_{i}f_{t,h}(\bold y)$ is in $C^1(\mathbb{R}^{n_{L+1}})$ and has
bounded first order derivatives.
Then, since $\bold{z}^{(L+1)}\,|\,\mathcal{F}_L$ is a centered Gaussian random vector with covariance matrix
\begin{equation}\label{eq:18052023primo}
\Sigma'_{n_{L+1}}:=(C_b+C_W\mathcal{O}_{\bold{n}_L}^{(L)})\mathrm{Id}_{n_{L+1}},
\end{equation}
by Lemma \ref{le:Gaussibp}$(ii)$ we have
\[
\mathbb E[z_i^{(L+1)}\partial_{i}f_{t,h}(\bold{z}^{(L+1)})\,|\,\mathcal{F}_L]=
\sum_{j=1}^{n_{L+1}}\Sigma'_{n_{L+1}}(i,j)\mathbb E[\partial_{ij}^2f_{t,h}(\bold{z}^{(L+1)})\,|\,\mathcal{F}_L].
\]
On combining this relation with \eqref{eq:09052023quarto} and taking the expectation,  we have
\begin{align}
\mathbb E [h_{t}(\bold{z}^{(L+1)})-h_{t}(\bold{z})]&=
\sum_{i,j=1}^{n_{L+1}}\mathbb{E}[(\Sigma'_{n_{L+1}}(i,j)-\Sigma_{n_{L+1}}(i,j))\mathbb E[\partial_{ij}^2f_{t,h}(\bold{z}^{(L+1)})\,|\,\mathcal{F}_L]]\nonumber\\
&=\sum_{i,j=1}^{n_{L+1}}\mathbb{E}[\mathbb E[(\Sigma'_{n_{L+1}}(i,j)-\Sigma_{n_{L+1}}(i,j))\partial_{ij}^2f_{t,h}(\bold{z}^{(L+1)})\,|\,\mathcal{F}_L]]\nonumber\\
&=\mathbb E\langle\Sigma'_{n_{L+1}}-\Sigma_{n_{L+1}},\mathrm{Hess}(f_{t,h}(\bold{z}^{(L+1)}))\rangle_{H.S.}\label{eq:09052023sesto}
\end{align}
where in the second equality we used the $\mathcal{F}_L$-measurability of $\mathcal{O}_{\bold{n}_L}^{(L)}$.
Taking the modulus on the relation \eqref{eq:09052023sesto} and applying the Cauchy-Schwarz inequality,  we have
\begin{align}
|\mathbb E [h_{t}(\bold{z}^{(L+1)})-h_{t}(\bold{z})]|
&\leq\sqrt{\mathbb E\|\Sigma'_{n_{L+1}}-\Sigma_{n_{L+1}}\|_{H.S.}^2}\sqrt{\mathbb E\|\mathrm{Hess}(f_{t,h}(\bold{z}^{(L+1)}))\|_{H.S.}^2}.
\label{eq:12052023primo}
\end{align}
By Lemma \ref{le:Stein2}$(ii)$,  for any $h\in\mathcal{I}_{n_{L+1}}$,  we have 
\begin{align*}
\mathbb E\mathbb \|\mathrm{Hess}(f_{t,h}(\bold{z}^{(L+1)}))\|_{H.S.}^2
\leq\|\Sigma_{n_{L+1}}^{-1}\|_{op}^2(n_{L+1}^2(\log t)^2\kappa
+530 n_{L+1}^{17/6}),\quad\text{$t\in (0,1)$.}
\end{align*}
Moreover,
\begin{equation}\label{eq:12052023secondo}
\mathbb E\|\Sigma'_{n_{L+1}}-\Sigma_{n_{L+1}}\|_{H.S.}^2=n_{L+1}C_W^2
\|\mathcal{O}_{\bold{n}_L}^{(L)}
-\mathcal{O}^{(L)}\|_{L^2}^2.
\end{equation}
On combining these relations and using the elementary inequality $\sqrt{a_1+a_2}\leq\sqrt a_1+\sqrt a_2$, $a_1,a_2\geq 0$, we have
\begin{align}
&\sup_{h\in\mathcal{I}_{n_{L+1}}}|\mathbb E [h_{t}(\bold{z}^{(L+1)})-h_{t}(\bold{z})]|\nonumber\\
&\qquad
\leq
C_W\|\Sigma_{n_{L+1}}^{-1}\|_{op}(n_{L+1}^{3/2}|\log t|\sqrt{\kappa
}+24 n_{L+1}^{23/12})\|\mathcal{O}_{\bold{n}_L}^{(L)}
-\mathcal{O}^{(L)}\|_{L^2}.\nonumber
\end{align}
On combining this latter inequality 
with \eqref{eq:02052023secondo},  we have
\begin{equation}\label{eq:02052023ottavo}
\kappa\leq\frac{4}{3}\|\Sigma_{n_{L+1}}^{-1}\|_{op}(n_{L+1}^{3/2}|\log t|\sqrt{\kappa}+24 n_{L+1}^{23/12})\gamma
+\frac{20 n_{L+1}}{\sqrt 2}\frac{\sqrt t}{1-t},\quad t\in (0,1).
\end{equation} 
Since $\kappa\leq 1$,  by this relation we have
\begin{equation*}
\kappa\leq\frac{4}{3}\|\Sigma_{n_{L+1}}^{-1}\|_{op}(n_{L+1}^{3/2}|\log t|+24 n_{L+1}^{23/12})\gamma
+\frac{20 n_{L+1}}{\sqrt 2}\frac{\sqrt t}{1-t},\quad t\in (0,1).
\end{equation*} 
Setting $t=\gamma^2$ in 
this latter inequality (note that this choice of the parameter $t$ is admissible since $\gamma\leq 1/\mathrm{e}<1$), we have
\begin{align}
\kappa&\leq\frac{4}{3}\|\Sigma_{n_{L+1}}^{-1}\|_{op}(2n_{L+1}^{3/2}|\log\gamma|+24 n_{L+1}^{23/12})\gamma
+\frac{20 n_{L+1}}{\sqrt 2}\frac{\gamma}{1-\gamma^2}\nonumber\\
&\leq\frac{4}{3}\|\Sigma_{n_{L+1}}^{-1}\|_{op}(2n_{L+1}^{3/2}|\log\gamma|+24 n_{L+1}^{23/12})\gamma
+20\sqrt 2 n_{L+1}\gamma,\label{eq:09052023settimo}
\end{align} 
where in the latter inequality we used the relation 
\begin{equation}\label{eq:02052023settimo}
\frac{1}{\sqrt 2}\frac{\gamma}{1-\gamma^2}\leq\sqrt 2\gamma,
\end{equation}
which holds since $\gamma\leq 1/\mathrm{e}<1/\sqrt 2$.
We rewrite the inequality \eqref{eq:09052023settimo} as
\begin{align*}
\kappa
&\leq\frac{8}{3}\|\Sigma_{n_{L+1}}^{-1}\|_{op}n_{L+1}^{3/2}\gamma|\log\gamma|+
(32 n_{L+1}^{23/12}\|\Sigma_{n_{L+1}}^{-1}\|_{op}
+20\sqrt 2 n_{L+1})\gamma.\nonumber
\end{align*} 
Taking the square root and multiplying by $|\log\gamma|$, we have
\begin{align*}
|\log\gamma|\sqrt{\kappa}
&\leq\sqrt{\frac{8}{3}}\|\Sigma_{n_{L+1}}^{-1}\|_{op}^{1/2}n_{L+1}^{3/4}\gamma^{1/2}|\log\gamma|^{3/2}+
(32 n_{L+1}^{23/12}\|\Sigma_{n_{L+1}}^{-1}\|_{op}
+20\sqrt 2 n_{L+1})^{1/2}\gamma^{1/2}|\log\gamma|.\nonumber
\end{align*} 
Since $\max\{\sup_{y\in (0,1/\mathrm e]}y^{1/2}|\log y|^{3/2},\sup_{y\in (0,1/\mathrm e]}y^{1/2}|\log y|^{1/2}\}\leq 4$,  we have
\begin{align}
|\log\gamma|\sqrt{\kappa}
&\leq 4\sqrt{\frac{8}{3}}\|\Sigma_{n_{L+1}}^{-1}\|_{op}^{1/2}n_{L+1}^{3/4}+
4(32 n_{L+1}^{23/12}\|\Sigma_{n_{L+1}}^{-1}\|_{op}
+20\sqrt 2 n_{L+1})^{1/2}\nonumber\\
&= 8\frac{\sqrt{6}}{3}\|\Sigma_{n_{L+1}}^{-1}\|_{op}^{1/2}n_{L+1}^{3/4}+
4(32 n_{L+1}^{23/12}\|\Sigma_{n_{L+1}}^{-1}\|_{op}
+20\sqrt 2 n_{L+1})^{1/2}\nonumber\\
&\leq 8\frac{\sqrt{6}}{3}\|\Sigma_{n_{L+1}}^{-1}\|_{op}^{1/2}n_{L+1}^{3/4}+
16\sqrt{2}n_{L+1}^{23/24}\|\Sigma_{n_{L+1}}^{-1}\|_{op}^{1/2}
+\sqrt{20\sqrt{2}}n_{L+1}^{1/2}\nonumber\\
&\leq\left[\left(8\frac{\sqrt{6}}{3}+16\sqrt{2}\right)\|\Sigma_{n_{L+1}}^{-1}\|_{op}^{1/2}
+\sqrt{20\sqrt{2}}\right]n_{L+1}^{23/24}\nonumber\\
&\leq (30\|\Sigma_{n_{L+1}}^{-1}\|_{op}^{1/2}
+6)n_{L+1}^{23/24}.\label{eq:02052023dieci}
\end{align}
By \eqref{eq:02052023ottavo} with $t=\gamma^2$,  \eqref{eq:02052023settimo} and \eqref{eq:02052023dieci},  we finally have
\eqref{eq:02052023undici},  indeed
\begin{align}
\kappa&\leq\frac{4}{3}\|\Sigma_{n_{L+1}}^{-1}\|_{op}(2n_{L+1}^{3/2}|\log\gamma|\sqrt{\kappa}+24 n_{L+1}^{23/12})\gamma
+20\sqrt{2}n_{L+1}\gamma\nonumber\\
&\leq\frac{4}{3}\|\Sigma_{n_{L+1}}^{-1}\|_{op}\left[\left(60\|\Sigma_{n_{L+1}}^{-1}\|_{op}^{1/2}+12\right)n_{L+1}^{59/24}+24n_{L+1}^{23/12}\right]\gamma
+20\sqrt{2}n_{L+1}\gamma\nonumber\\
&\leq(80\|\Sigma_{n_{L+1}}^{-1}\|_{op}^{3/2}+48\|\Sigma_{n_{L+1}}^{-1}\|_{op}+20\sqrt{2})n_{L+1}^{59/24}\gamma.\nonumber
\end{align}
The proof is completed.
\\
\noindent$\square$

\subsection{Normal approximation of deep random Gaussian NNs in the $1$-Wasserstein distance}

The following theorem holds.

\begin{Theorem}\label{thm:1Wasserstein}
Let $\mathrm{GNN}(L,n_0,n_{L+1},\bold{n}_{L},\sigma,\bold x,\bold b,\bold W)$ be a deep random Gaussian NN,  and let
the notation and assumptions of Theorem \ref{le:fund_est} prevail.  Then
\begin{align}
d_{W_1}(\bold{z}^{(L+1)},\bold{z})&\leq K_1\left(\sum_{k=1}^{L}[4\sqrt{2}\|P(|Z|)\|_{L^2}]^{L-k}\frac{c_k}{\sqrt{n_k}}\right)\nonumber\\
&\leq K_i\left(\sum_{k=1}^{L}[4\sqrt{2}\|P(|Z|)\|_{L^2}]^{L-k}\frac{c_k}{\sqrt{n_k}}\right),\quad\text{$i=2,3$}\nonumber
\end{align}
where the constants $c_k$, $k=1,\ldots,L$,  are defined by \eqref{eq:c20062023},
\[
K_1=K_1(n_0,n_{L+1},\bold x, \sigma,C_b,C_W):=\frac{n_{L+1}C_W}{\sqrt{C_b+C_W\mathcal{O}^{(L)}}},
\]
\[
K_2=K_2(n_{L+1},C_b,C_W):=\frac{n_{L+1}C_W}{\sqrt{C_b}}
\]
and
\[
K_3=K_3(n_0,n_{L+1},\bold x, \sigma,C_W):=\frac{n_{L+1}\sqrt{C_W}}{\sqrt{\mathcal{O}^{(L)}}},
\]
\end{Theorem}

\begin{Remark}\label{re:19092023} 
Let $\mathrm{GNN}(L,n_0,1,\bold{n}_L,\sigma,\bold x,\bold b,\bold W)$ be
a deep random Gassian NN with univariate output,  widths of the hidden layers satsfying \eqref{eq:19092023} and a polynomially bounded to order $r\geq 1$ activation function $\sigma$.  Then by Theorem 3.3 in \cite{FHMNP} we have that there exist two constants $C,C_0>0$ such that
\[
\frac{C_0}{n}\leq d_{W_1}(z^{(L+1)},z)\leq\frac{C}{n}.
\]
Clearly this inequality shows the optimality of the rate $1/n$.  Here again,  we note that the corresponding bound provided by Theorem \ref{thm:1Wasserstein} gives (under different assumptions on $\sigma$ and without any condition on the widths of the hidden layers) a detailed description of the analytical dependence of the upper estimate on the parameters of the model.  This makes our result useful for the purpose of output localization (see Section \ref{sec:illustrations}). 
\end{Remark}

\noindent{\it Proof.}
Let $g\in\mathcal{L}_{n_{L+1}}(1)$ be arbitrarily fixed.  Without loss of generality we assume that $\bold{z}$ is independent of $\mathcal{F}_L$.
By Lemma \ref{le:Stein} we then have
\begin{align}
\mathbb E [g(\bold z^{(L+1)})-g(\bold z)\,|\,\mathcal{F}_L]
&=\sum_{i=1}^{n_{L+1}}\mathbb E[z_i^{(L+1)}\partial_{i}f_{g}(\bold{z}^{(L+1)})\,|\,\mathcal{F}_L]
-\sum_{i,j=1}^{n_{L+1}}\Sigma_{n_{L+1}}(i,j)\mathbb{E}[\partial_{ij}^2 f_{g}(\bold{z}^{(L+1)})\,|\,\mathcal{F}_L],
\label{eq:09052023dodici}
\end{align}
where
\[
f_{g}(\bold y):=\int_0^\infty\mathbb{E}[g(\bold z)-g(\mathrm{e}^{-t}\bold y+\sqrt{1-\mathrm{e}^{-2t}}\bold z)]\,\dd t,\quad \bold y\in\mathbb R^{n_{L+1}}.
\]
Again by Lemma \ref{le:Stein} we have that
the mapping $\bold y\mapsto\partial_{i}f_{g}(\bold y)$ is in $C^1(\mathbb{R}^{n_{L+1}})$ and has
bounded first order derivatives. 
Applying Lemma \ref{le:Gaussibp}$(ii)$ exactly as in the proof of Theorem \ref{thm:convexdistance} (see a few lines before Equation \eqref{eq:09052023sesto}),  we have
\[
\mathbb E [g(\bold{z}^{(L+1)})-g(\bold{z})]
=\mathbb E\langle\Sigma'_{n_{L+1}}-\Sigma_{n_{L+1}},\mathrm{Hess}(f_{g}(\bold{z}^{(L+1)}))\rangle_{H.S.},
\]
where the matrix $\Sigma'_{n_{L+1}}$ is defined by \eqref{eq:18052023primo}.
Therefore,  applying the Cauchy-Schwarz inequality as in \eqref{eq:12052023primo},
we have
\begin{align}
|\mathbb E [g(\bold{z}^{(L+1)})-g(\bold{z})]|
&\leq\sqrt{\mathbb E\|\Sigma'_{n_{L+1}}-\Sigma_{n_{L+1}}\|_{H.S.}^2}\sqrt{\mathbb E\|\mathrm{Hess}(f_{g}(\bold{z}^{(L+1)}))\|_{H.S.}^2}
\end{align}
By Lemma \ref{le:Stein} we have
\[
\sqrt{\mathbb E\|\mathrm{Hess}(f_{g}(\bold{z}^{(L+1)}))\|_{H.S.}^2}\leq\sup_{\bold{y}\in\mathbb R^{n_{L+1}}}\|\mathrm{Hess}f_{g}(\bold y)\|_{H.S.}\leq\sqrt{n_{L+1}}\|\Sigma_{n_{L+1}}^{-1}\|_{op}\|\Sigma_{n_{L+1}}\|_{op}^{1/2},
\]
On combining these relations with \eqref{eq:12052023secondo},  we have
\begin{align*}
|\mathbb E [g(\bold{z}^{(L+1)})-g(\bold{z})]|
&\leq n_{L+1}C_W\|\Sigma_{n_{L+1}}^{-1}\|_{op}\|\Sigma_{n_{L+1}}\|_{op}^{1/2}\|\mathcal{O}_{\bold{n}_L}^{(L)}-\mathcal{O}^{(L)}\|_{L^2}.
\end{align*}
The claim follows taking the supremum over $g$ on this inequality and then using Theorem \ref{le:fund_est},  relation \eqref{eq:10052023secondo} and the fact that
\[
\|\Sigma_{n_{L+1}}\|_{op}=C_b+C_W\mathcal{O}^{(L)}.
\]
\noindent$\square$

\section{Localization of the output}\label{sec:illustrations}
In this section we want to show the potentiality of the obtained results for practical applications. Indeed, 
both Theorem~\ref{prop:shallow}$(ii)$ and Theorem~\ref{thm:convexdistance}
allow one to explicitly estimate the probability that the output of a random
Gaussian NN evaluated at the input $\mathbf x$ belongs to a suitable set,
without resorting to computationally expensive Monte Carlo methods. 
This is what we call output localization,  which can suggest a suitable architecture design to estimate a target function $f$.
Indeed in statistical learning, 
given a training set
\[
\{ (\mathbf x_i, f(\mathbf x_i )\}_{i\in I}\subset \mathbb R^{n_0}\times \mathbb R^{n_{L+1}},
\]
the goal is to estimate the unknown function $f$ by a NN with a certain fixed architecture. This is performed
by minimizing an empirical risk function over the space of NN's parameters, i.e.,  the biases and the weights.  Such kind of minimization problems are non-convex and are usually
addressed by a gradient descent procedure; the latter stabilizes towards a solution which depends on the inizialization point.  Since in practical applications the inizialization point is given by a realization of a random Gaussian NN,
it can be convenient to choose $L$, $\bold{n}_L$, $\sigma$, $C_b$ and $C_W$
by exploiting output localization,  i.e.,  in such a way to have a good estimate of
the probability
$\mathbb P(\mathbf z^{(L+1)}\left(\mathbf x_i\right)\in V_i)$,  $i\in I$,  where 
$V_i$ is an appropriate neighborhood of $f(\mathbf x_{i})$.

Hence, in the following we want to show how to use the upper bound on the total variation distance provided in Theorem \ref{prop:shallow}$(ii)$,  for a shallow random Gaussian NN with univariate output,  and the upper bound on the convex distance provided by Theorem \ref{thm:convexdistance},  for a deep random Gaussian NN,  to localize the output.  Indeed,  given a measurable convex set $V\subset\mathbb{R}^{n_{L+1}}$,  by the definitions of both the total variation and the convex distances we have
\begin{equation*}
\mathbb{P}(\bold{z}\in V)-C_{bound} \leq\mathbb{P}(\bold{z}^{(L+1)}\in V)\leq\mathbb{P}(\bold{z}\in V)+C_{bound},
\end{equation*}
where $L=1$, $n_2=1$ and
\begin{equation} 
\label{eq:C_{TVbound}}
C_{bound}:=2\frac{C_W\sqrt{\mathbb{V}\mathrm{ar}(\sigma(Z\sqrt{C_b+C_W\mathcal{O}^{(0)}})^2)}}{C_b+C_W\mathbb{E}\sigma(Z\sqrt{C_b+C_W\mathcal{O}^{(0)}})^2}\frac{1}{\sqrt n_1},
\end{equation}
for the case of a shallow NN,  while
\begin{equation} 
\label{eq:C_{bound}}
C_{bound}:=C_1 \left(\sum_{k=1}^{L}[4\sqrt{2}\|P(|Z|)\|_{L^2}]^{L-k}\frac{c_k}{\sqrt{n_k}}\right),
\end{equation}
for the case of a deep NN,  where $C_1$ is given by \eqref{eq: C_1 esplicita} and the constants $c_k$, $k=1,\ldots,L$,  are defined by \eqref{eq:c20062023}.

Let $V:=\prod_{i=1}^{n_{L+1}}[r_i,s_i]$ be a rectangle of $\mathbb R^{n_{L+1}}$.
Since $\bold{z}$ is an $n_{L+1}$-dimensional centered Gaussian random vector with covariance matrix \eqref{eq:09052023primo},  we have
\[
\mathbb{P}(\bold{z}\in V)=\prod_{i=1}^{n_{L+1}}\left(\mathbb{P}\left(Z\leq\frac{s_i}{\sqrt{C_b+C_W\mathcal{O}^{(L)}}}\right)-\mathbb{P}\left(Z\leq\frac{r_i}{\sqrt{C_b+C_W\mathcal{O}^{(L)}}}\right)\right).
\]
\begin{table}[!h]
\centering
\small
\resizebox{\textwidth}{!}{
\begin{tabular}{l|rrr|rrr|rrr|rrr|rrr|rrr}
\hline
\hline
\multicolumn{19}{c}{$\mathbf{ x=(0,0,0,0)}$}\\
\\
\hline
$n$ & \multicolumn{3}{c}{$1$}&\multicolumn{3}{c}{$10$}&\multicolumn{3}{c}{$10^2$}&\multicolumn{3}{c}{$10^3$}&\multicolumn{3}{c}{$10^4$}&\multicolumn{3}{c}{$10^5$}\\
  \hline
\diagbox{$C_{b}$}{$C_{W}$} & 0.01 & 0.1 & 1 & 0.01 & 0.1 & 1 & 0.01 & 0.1 & 1 & 0.01 & 0.1 & 1 & 0.01 & 0.1 & 1& 0.01 & 0.1 & 1 \\ 
  \hline
1 & 0.02 & 0.21 & 1.49 & 0.01 & 0.07 & 0.47 & 0.00 & 0.02 & 0.15 & 0.00 & 0.01 & 0.05 & 0.00 & 0.00 & 0.01 & 0.00 & 0.00 & 0.00 \\ 
  10 & 0.01 & 0.07 & 0.47 & 0.00 & 0.02 & 0.15 & 0.00 & 0.01 & 0.05 & 0.00 & 0.00 & 0.01 & 0.00 & 0.00 & 0.00 & 0.00 & 0.00 & 0.00 \\   \hline
   \hline
\multicolumn{19}{c}{$\mathbf{x=(0.1,0.1,0.1,0.1)}$}\\
\\
\hline
$n$ & \multicolumn{3}{c}{$1$}&\multicolumn{3}{c}{$10$}&\multicolumn{3}{c}{$10^2$}&\multicolumn{3}{c}{$10^3$}&\multicolumn{3}{c}{$10^4$}&\multicolumn{3}{c}{$10^5$}\\
  \hline
\diagbox{$C_{b}$}{$C_{W}$} & 0.01 & 0.1 & 1 & 0.01 & 0.1 & 1 & 0.01 & 0.1 & 1 & 0.01 & 0.1 & 1 & 0.01 & 0.1 & 1& 0.01 & 0.1 & 1 \\ 
  \hline
1 &  0.02 & 0.21 & 1.49 & 0.01 & 0.07 & 0.47 & 0.00 & 0.02 & 0.15 & 0.00 & 0.01 & 0.05 & 0.00 & 0.00 & 0.01 & 0.00 & 0.00 & 0.00 \\ 
  10 & 0.01 & 0.07 & 0.47 & 0.00 & 0.02 & 0.15 & 0.00 & 0.01 & 0.05 & 0.00 & 0.00 & 0.01 & 0.00 & 0.00 & 0.00 & 0.00 & 0.00 & 0.00 \\ 
   \hline
   \hline
\multicolumn{19}{c}{$\mathbf{x=(0.5,-0.5,0.5,-0.5)}$}\\
\\
\hline
$n$ &\multicolumn{3}{c}{$1$}&\multicolumn{3}{c}{$10$}&\multicolumn{3}{c}{$10^2$}&\multicolumn{3}{c}{$10^3$}&\multicolumn{3}{c}{$10^4$}&\multicolumn{3}{c}{$10^5$}\\
  \hline
\diagbox{$C_{b}$}{$C_{W}$} & 0.01 & 0.1 & 1 & 0.01 & 0.1 & 1 & 0.01 & 0.1 & 1 & 0.01 & 0.1 & 1 & 0.01 & 0.1 & 1& 0.01 & 0.1 & 1 \\ 
  \hline
1 & 0.02 & 0.22 & 1.54 & 0.01 & 0.07 & 0.49 & 0.00 & 0.02 & 0.15 & 0.00 & 0.01 & 0.05 & 0.00 & 0.00 & 0.02 & 0.00 & 0.00 & 0.00 \\ 
  10 & 0.01 & 0.07 & 0.47 & 0.00 & 0.02 & 0.15 & 0.00 & 0.01 & 0.05 & 0.00 & 0.00 & 0.01 & 0.00 & 0.00 & 0.00 & 0.00 & 0.00 & 0.00 \\
   \hline
   \hline
\multicolumn{19}{c}{$\mathbf{x=(10,10,10,10)}$}\\
\\
\hline
$n$ & \multicolumn{3}{c}{$1$}&\multicolumn{3}{c}{$10$}&\multicolumn{3}{c}{$10^2$}&\multicolumn{3}{c}{$10^3$}&\multicolumn{3}{c}{$10^4$}&\multicolumn{3}{c}{$10^5$}\\
  \hline
\diagbox{$C_{b}$}{$C_{W}$} & 0.01 & 0.1 & 1 & 0.01 & 0.1 & 1 & 0.01 & 0.1 & 1 & 0.01 & 0.1 & 1 & 0.01 & 0.1 & 1& 0.01 & 0.1 & 1 \\ 
  \hline
1 & 0.03 & 0.48 & 0.44 & 0.01 & 0.15 & 0.14 & 0.00 & 0.05 & 0.04 & 0.00 & 0.02 & 0.01 & 0.00 & 0.00 & 0.00 & 0.00 & 0.00 & 0.00 \\ 
  10 & 0.01 & 0.09 & 0.36 & 0.00 & 0.03 & 0.11 & 0.00 & 0.01 & 0.04 & 0.00 & 0.00 & 0.01 & 0.00 & 0.00 & 0.00 & 0.00 & 0.00 & 0.00 \\
   \hline
   \hline
\end{tabular}
}
\caption{Values of $C_{bound}$ in \eqref{eq:C_{TVbound}}
for different inputs $\bold x$ for a shallow random Gaussian NN with architecture $L=1$, $n_0=4$, $n_1=n$, $n_2=1$ and ReLu activation function.}
 \label{tab:C-TVbound}
\end{table}

\begin{table}[!h] 
\centering
\small
\resizebox{\textwidth}{!}{
\begin{tabular}{l|rrr|rrr|rrr|rrr|rrr|rrr}
\hline
\hline
\multicolumn{19}{c}{$\mathbf{ x=(0,0,0,0)}$}\\
\\
\hline
$n$ & \multicolumn{3}{c}{$10^4$}&\multicolumn{3}{c}{$10^5$}&\multicolumn{3}{c}{$10^6$}&\multicolumn{3}{c}{$10^7$}&\multicolumn{3}{c}{$10^8$}&\multicolumn{3}{c}{$10^9$}\\
  \hline
\diagbox{$C_{b}$}{$C_{W}$} & 0.01 & 0.1 & 1 & 0.01 & 0.1 & 1 & 0.01 & 0.1 & 1 & 0.01 & 0.1 & 1 & 0.01 & 0.1 & 1& 0.01 & 0.1 & 1 \\ 
  \hline
1 & 0.03 & 0.58 & 88.59 & 0.01 & 0.18 & 28.01 & 0.00 & 0.06 & 8.86 & 0.00 & 0.02 & 2.80 & 0.00 & 0.01 & 0.89 & 0.00 & 0.00 & 0.28 \\ 
  10 & 0.07 & 1.38 & 331.57 & 0.02 & 0.44 & 104.85 & 0.01 & 0.14 & 33.16 & 0.00 & 0.04 & 10.49 & 0.00 & 0.01 & 3.32 & 0.00 & 0.00 & 1.05 \\    \hline
   \hline
\multicolumn{19}{c}{$\mathbf{x=(0.1,0.1,0.1,0.1)}$}\\
\\
\hline
$n$ & \multicolumn{3}{c}{$10^4$}&\multicolumn{3}{c}{$10^5$}&\multicolumn{3}{c}{$10^6$}&\multicolumn{3}{c}{$10^7$}&\multicolumn{3}{c}{$10^8$}&\multicolumn{3}{c}{$10^9$}\\
  \hline
\diagbox{$C_{b}$}{$C_{W}$} & 0.01 & 0.1 & 1 & 0.01 & 0.1 & 1 & 0.01 & 0.1 & 1 & 0.01 & 0.1 & 1 & 0.01 & 0.1 & 1& 0.01 & 0.1 & 1 \\ 
  \hline
 1 & 0.03 & 0.58 & 89.30 & 0.01 & 0.18 & 28.24 & 0.00 & 0.06 & 8.93 & 0.00 & 0.02 & 2.82 & 0.00 & 0.01 & 0.89 & 0.00 & 0.00 & 0.28 \\ 
  10 & 0.07 & 1.38 & 331.85 & 0.02 & 0.44 & 104.94 & 0.01 & 0.14 & 33.18 & 0.00 & 0.04 & 10.49 & 0.00 & 0.01 & 3.32 & 0.00 & 0.00 & 1.05 \\ 
   \hline
   \hline
\multicolumn{19}{c}{$\mathbf{x=(0.5,-0.5,0.5,-0.5)}$}\\
\\
\hline
$n$ &\multicolumn{3}{c}{$10^4$}&\multicolumn{3}{c}{$10^5$}&\multicolumn{3}{c}{$10^6$}&\multicolumn{3}{c}{$10^7$}&\multicolumn{3}{c}{$10^8$}&\multicolumn{3}{c}{$10^9$}\\
  \hline
\diagbox{$C_{b}$}{$C_{W}$} & 0.01 & 0.1 & 1 & 0.01 & 0.1 & 1 & 0.01 & 0.1 & 1 & 0.01 & 0.1 & 1 & 0.01 & 0.1 & 1& 0.01 & 0.1 & 1 \\ 
  \hline
 1 & 0.03 & 0.58 & 106.14 & 0.01 & 0.18 & 33.56 & 0.00 & 0.06 & 10.61 & 0.00 & 0.02 & 3.36 & 0.00 & 0.01 & 1.06 & 0.00 & 0.00 & 0.34 \\ 
  10 & 0.07 & 1.38 & 338.62 & 0.02 & 0.44 & 107.08 & 0.01 & 0.14 & 33.86 & 0.00 & 0.04 & 10.71 & 0.00 & 0.01 & 3.39 & 0.00 & 0.00 & 1.07 \\ 
   \hline
   \hline
\multicolumn{19}{c}{$\mathbf{x=(10,10,10,10)}$}\\
\\
\hline
$n$ & \multicolumn{3}{c}{$10^4$}&\multicolumn{3}{c}{$10^5$}&\multicolumn{3}{c}{$10^6$}&\multicolumn{3}{c}{$10^7$}&\multicolumn{3}{c}{$10^8$}&\multicolumn{3}{c}{$10^9$}\\
  \hline
\diagbox{$C_{b}$}{$C_{W}$} & 0.01 & 0.1 & 1 & 0.01 & 0.1 & 1 & 0.01 & 0.1 & 1 & 0.01 & 0.1 & 1 & 0.01 & 0.1 & 1& 0.01 & 0.1 & 1 \\ 
  \hline
 1 & 0.03 & 1.90 & 2998.50 & 0.01 & 0.60 & 948.21 & 0.00 & 0.19 & 299.85 & 0.00 & 0.06 & 94.82 & 0.00 & 0.02 & 29.99 & 0.00 & 0.01 & 9.48 \\ 
  10 & 0.07 & 1.69 & 3027.47 & 0.02 & 0.54 & 957.37 & 0.01 & 0.17 & 302.75 & 0.00 & 0.05 & 95.74 & 0.00 & 0.02 & 30.27 & 0.00 & 0.01 & 9.57 \\ 
   \hline
\end{tabular}
}
\caption{Values of $C_{bound}$ in \eqref{eq:C_{bound}} for different inputs $\bold x$ for a deep random Gaussian NN with architecture $L=3$, $n_0=4$, $n_1=n_2=n_3=n$, $n_4=1$ and ReLu activation function.}
 \label{tab:C-bound}
\end{table}

\vspace{1cm}

\begin{table}[ht] 
\centering
\small
\resizebox{\textwidth}{!}{
\begin{tabular}{l|rrr|rrr|rrr|rrr}
 \hline
& \multicolumn{3}{c}{$\mathbf{x=(0,0,0,0)}$}&\multicolumn{3}{c}{$\mathbf{x=(0.1,0.1,0.1,0.1)}$}&\multicolumn{3}{c}{$\mathbf{x=(0.5,-0.5,0.5,-0.5)}$}&\multicolumn{3}{c}{$\mathbf{x=(10,10,10,10)}$}\\
  \hline
\diagbox{$C_{b}$}{$C_{W}$} & 0.01 & 0.1 & 1 & 0.01 & 0.1 & 1 & 0.01 & 0.1 & 1 & 0.01 & 0.1 & 1 \\ 
  \hline
1 & 1.55 & 14.80 & 85.04 & 1.55 & 14.80 & 85.00 & 1.55 & 14.80 & 83.86 & 1.55 & 14.78 & 33.09 \\ 
  10 & 0.36 & 3.52 & 31.83 & 0.36 & 3.52 & 31.83 & 0.36 & 3.52 & 31.82 & 0.36 & 3.52 & 30.28 \\ 
   \hline
\end{tabular}
}
\caption{Values of $C_1$ in \eqref{eq: C_1 esplicita} for different inputs $\bold x$ 
for a deep random Gaussian NN with architecture $L=3$, $n_0=4$, $n_1=n_2=n_3=n$, $n_4=1$ and ReLu activation function.
}
\label{tab:C1} 
\end{table}
%
Now,  we furnish numerical values of the constant $C_{bound}$ in \eqref{eq:C_{TVbound}} and \eqref{eq:C_{bound}}; in both cases we take $\sigma(x):=x\bold{1}\{x\geq 0\}$, i.e.,  a ReLu activation function.
Since the ReLu function is Lipschitz continuous with Lipschitz constant equal to $1$,  $\mathbb E Z^2=1$ and $\mathbb E Z^4=3$,  by Proposition \ref{prop:04052023primo}$(ii)$ we have
\[
\|P(|Z|)\|_{L^2}=C_W\sqrt{\mathbb E Z^4\bold{1}\{Z\geq 0\}}=C_W\sqrt{3/2}.
\]
By the expression of the constants $c_\ell$ in Theorem \ref{le:fund_est},  we have
\[
c_\ell=(C_b+C_W\mathcal{O}^{(\ell-1)})\sqrt{2\mathbb E Z^4\bold{1}\{Z\geq 0\}}
=(C_b+C_W\mathcal{O}^{(\ell-1)})\sqrt{3},\quad\text{$\ell=1,\ldots,L$.}
\]
By the definition of the quantities $\mathcal{O}^{(\ell)}$,  $\ell=1,\ldots,L$,  (see \eqref{eq:20062023sesto}) we have
\[
\mathcal{O}^{(\ell)}=(C_b+C_W\mathcal{O}^{(\ell-1)})\mathbb{E}Z^2\bold{1}\{Z\geq 0\}
=(C_b+C_W\mathcal{O}^{(\ell-1)})/2,\quad\text{$\ell=1,\ldots,L$}
\]
with $\mathcal{O}^{(0)}$ given by \eqref{eq:20062023quinto}.
Therefore,
\[
\mathcal{O}^{(\ell)}=\frac{C_b}{2}\sum_{k=0}^{\ell-1}\frac{C_W^k}{2^k}+\frac{C_W^{\ell}}{2^{\ell}}\mathcal{O}^{(0)}.
\]
Table \ref{tab:C-TVbound}  gives the values of $C_{bound}$ in \eqref{eq:C_{TVbound}} in the case of a shallow random Gaussian NN with architecture $L=1$,  $n_0=4$,  $n_1=n$ and $n_2=1$, 
for different values of $n\in\{1, 10,10^2,10^3,10^4,10^5  \}$.
Table \ref{tab:C-bound} gives the values of $C_{bound}$ in \eqref{eq:C_{bound}} in the case of a deep random Gaussian NN with architecture $L=3$, $n_0=4$,  $n_1=n_2=n_3=n$ and $n_4=1$,  for different values of $n \in \{10^4, 10^5,10^6,10^7,10^8,10^9  \}$.  In both tables we consider four different inputs
\[
\bold x\in\{(0, 0, 0, 0), (0.1, 0.1, 0.1, 0.1),(0.5, -0.5, 0.5, -0.5),(10, 10, 10, 10)\},
\]
whose Euclidean norm is, respectively,  $0$,  strictly less than $1$,  $1$ and strictly larger than $1$,  and
\[
C_b\in\{ 1, 10 \}\quad\text{and}\quad C_W\in\{  0.01 , 0.1, 1 \}.
\]
It is clear from Table \ref{tab:C-TVbound} that for a shallow random Gaussian NN,  the Gaussian approximation of the output is very good for any of the choices of the parameters $n$,  $C_b$,  $C_W$ and of the input $\bold x$.  It is also clear from Table \ref{tab:C-bound} that
for a deep random Gaussian NN,  such as the considered NN with three hidden layers,  
the Gaussian approximation of the output is very good for some choices of the parameters $C_b$,  $C_W$,  $n$ and of the input $\bold x$ (also when the number $n$ of neurons in the hidden layers is not excessive) but it is very poor for other choices of these quantities. 
In order to analyse the influence of the parameters $C_b$,  $C_W$ and of the input $\bold x$ on the value of $C_{bound}$ in \eqref{eq:C_{bound}},
in Table \ref{tab:C1} we reported the value of the constant $C_1$ that appears in \eqref{eq:C_{bound}} and it is explicitly given in \eqref{eq: C_1 esplicita}.  
As it can be seen from this table,  the value of $C_{bound}$ strongly depends on the value of $C_1$ which,  in turn,  it is closely related to the choice of the parameters $C_b$ and $C_W$ and does not depend much on the norm of the input vector $\bold x$.

\end{document}